\documentclass[12pt,a4paper,dvipdfmx]{amsart}

\usepackage{amsmath,amsfonts,amsthm,amssymb,amscd,extarrows}

\usepackage{graphicx}
\usepackage{here}   
\usepackage[stable]{footmisc}
\usepackage{cancel}
\usepackage{ascmac,fancybox}
\usepackage{cases}   
\usepackage{ulem}   
\usepackage[all]{xy}
\usepackage{url}
\usepackage{color}

\makeatletter
  
  \@addtoreset{equation}{section}
\makeatother
\newtheorem{Theorem}{Theorem}[section]
\newtheorem{Proposition}[Theorem]{Proposition}
\newtheorem{Lemma}[Theorem]{Lemma}
\newtheorem{Corollary}[Theorem]{Corollary}
\newtheorem{Remark}[Theorem]{Remark}

\newtheorem{Definition}[Theorem]{Definition}

\newtheorem{Claim}[Theorem]{Claim}

\def\RMN#1{\uppercase\expandafter{\romannumeral#1}}

\newcommand{\spec}{\mathop{\rm Spec}\nolimits}
\newcommand{\proj}{\mathop{\rm Proj}\nolimits}

\newcommand{\bfp}{\mathfrak{p}}

\newcommand{\lcal}{\mathcal{L}}
\newcommand{\oo}{\mathcal{O}}

\newcommand{\bZ}{\mathbb{Z}}
\newcommand{\bR}{\mathbb{R}}
\newcommand{\bRo}{\mathbb{R}_{\ge 0}}

\newcommand{\bP}{\mathbb{P}}
\newcommand{\bN}{\mathbb{N}}
\newcommand{\bNo}{\mathbb{N}_0}

\newcommand{\Proof}{\noindent {\it Proof.} \ }

\newcommand{\niretu}[2]{
\tiny \begin{array}{c} {#1} \\ {#2} \end{array}
}

\newcommand{\sanretu}[3]{
\tiny \begin{array}{c} {#1} \\ {#2} \\ {#3} \end{array}
}

\begin{document}
\title[finite generation 
of symbolic Rees rings]{
Some necessary and sufficient condition for finite generation 
of symbolic Rees rings}
\author{Taro Inagawa and Kazuhiko Kurano}
\dedicatory{}
\date{}
\thanks{The first author was supported by 
 JSPS KAKENHI Grant Number 19H00637. \\
The second author was supported by 
 JSPS KAKENHI Grant Numbers 18K03226, 19H00637. \\
 This work was partly supported by Osaka City University Advanced Mathematical Institute: MEXT Joint Usage/Research Center on Mathematics and Theoretical Physics JPMXP0619217849. }
\maketitle

%13A30, 14E99

\begin{abstract}
Consider the blow-up $Y$ of a weighted projective plane at a point in the open orbit over a field of characteristic $0$.
We assume that there exists a curve $C$ on $Y$ such that $C^2<0$ and $C.E=1$, where $E$ is the exceptional curve.

In this paper we give a (very simple) necessary and sufficient condition for finite generation of the Cox ring of $Y$ (Theorem~\ref{Keuprop}).
It is an affirmative answer to a conjecture due to He and Kurano-Nishida.
\end{abstract}

\section{Introduction}

For pairwise coprime positive integers
$a$, $b$ and $c$, 
let ${\frak p}$ be the defining ideal of the space monomial
curve $(t^a, t^b, t^c)$ in $K^3$, where $K$ is a field.
The ideal ${\frak p}$ is generated by at most three binomials in $P=K[x,y,z]$
(Herzog~\cite{Her}).
The symbolic Rees rings of space monomial primes are deeply studied by many authors. 
Huneke~\cite{Hu} and Cutkosky~\cite{C} developed criteria 
for finite generation of such rings. 
In 1994,
Goto-Nishida-Watanabe~\cite{GNW} first found examples of infinitely generated
symbolic Rees rings of space monomial primes.
Recently, using toric geometry,
 Gonz\'alez-Karu~\cite{GK} found some sufficient conditions for 
the symbolic Rees rings of space monomial primes to be infinitely generated.

Cutkosky~\cite{C} found the geometric meaning of the symbolic Rees rings of space monomial primes. 
Let ${\Bbb P}(a,b,c)$ be the weighted projective surface with degree $a$, $b$, $c$.
Let $Y$ be the blow-up at a point in the open orbit of the toric variety ${\Bbb P}(a,b,c)$.
Then the Cox ring of $Y$ is isomorphic to the extended symbolic Rees ring
of the space monomial prime ${\frak p}$.
Therefore, the symbolic Rees ring of the space monomial prime ${\frak p}$
is finitely generated if and only if the Cox ring of $Y$ is finitely generated,
that is, $Y$ is a Mori dream space.
A curve $C$ on $Y$ is called a negative curve if $C^2 < 0$ and
$C$ is different from the exceptional curve $E$. 
Here suppose $\sqrt{abc} \not\in {\Bbb Q}$.
Cutkosky~\cite{C} proved that the symbolic Rees ring of the space monomial prime ${\frak p}$ is finitely generated if and only if
the following two conditions are satisfied:
\begin{itemize}
\item[(1)]
There exists a negative curve $C$.
\item[(2)]
There exists a curve $D$ on $Y$ such that $C \cap D =\emptyset$.
\end{itemize}
In the case of ${\rm ch}(K) >0$, Cutkosky~\cite{C} proved that 
the symbolic Rees ring is Noetherian if there exists a negative curve.

The existence of negative curves is a very difficult and important problem,
that is deeply related to the Nagata conjecture (Proposition~5.2 in Cutkosky-Kurano~\cite{CK}) and the rationality of Seshadri constant.
The existence of negative curves is studied in 
Gonz\'alez-AnayaGonz\'alez-Karu~\cite{GAGK2}, \cite{GAGK3}, Kurano-Matsuoka~\cite{KM} and Kurano~\cite{K42}.

Examples that have a negative curve $C$ such that $C.E \ge 2$ 
are studied in 
Gonz\'alez-AnayaGonz\'alez-Karu~\cite{GAGK}, \cite{GAGK1} and Kurano-Nishida~\cite{KN}.

In this paper, we study the case where 
there exists a negative curve $C$ such that $C.E = 1$ in the case of ${\rm ch}(K) =0$.

Now, we state the main result of this paper precisely.

Let $K$ be a field.
Suppose that $a$, $b$, $c$ are pairwise coprime positive integers.
Let ${\frak p}$ be the kernel of the $K$-algebra map
$\varphi : P = K[x,y,z] \rightarrow K[T]$ defined by 
$\varphi(x) = T^a$, $\varphi(y) = T^b$, $\varphi(z) = T^c$.
Assume that ${\frak p}$ is not complete intersection.
Then we know
\[
{\frak p} 
%= I_2\left(
%\begin{array}{ccc}
%x^{s_2} & y^{t_3} & z^{u_1} \\
%y^{t_1} & z^{u_2} & x^{s_3}
%\end{array}
%\right) 
= (x^s-y^{t_1}z^{u_1}, y^t-x^{s_2}z^{u_2}, z^u-x^{s_3}y^{t_3})
\]
with positive integers $s_2$, $s_3$, $t_1$, $t_3$, $u_1$, $u_2$
such that $s = s_2+s_3$, $t = t_1+t_3$, $u = u_1+u_2$.
One can prove ${\rm GCD}(s_2,s_3) = {\rm GCD}(t_1,t_3) = {\rm GCD}(u_1,u_2) = 1$ as in the proof of Proposition~4.8 in \cite{KN}.
We put $\overline{t} = -t/t_3$, $\overline{u}=-u_2/u$, $\overline{s}=s_2/s_3$.
Remark $\overline{t} < -1 < \overline{u} < 0 < \overline{s}$.
Here consider the triangle $\Delta_{\overline{t}, \overline{u}, \overline{s}}$ as follows:
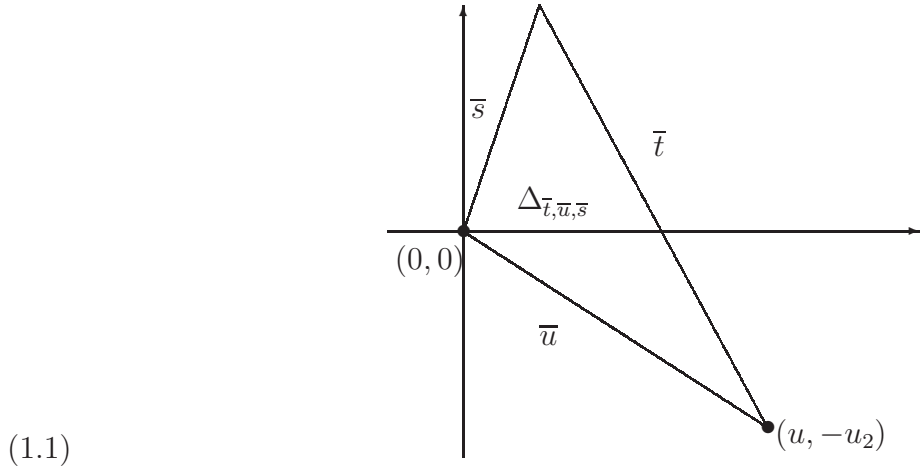
\begin{equation}\label{tri}
{
\setlength\unitlength{1truecm}
  \begin{picture}(7,6)(-2,-3)
  \put(-1,0){\vector(1,0){7}}
  \put(0,-3){\vector(0,1){6}}
\qbezier (0,0) (2,-1.3) (4,-2.6)
\qbezier (0,0) (0.5,1.5) (1,3)
\qbezier (1,3) (2.5,0.2) (4,-2.6)
\put(-0.1,-0.1){$\bullet$}
\put(3.9,-2.7){$\bullet$}
\put(-0.9,-0.5){$(0,0)$}
\put(4.1,-2.8){$(u, -u_2)$}
  \put(1,-1.5){$\overline{u}$}
  \put(2.5,1){$\overline{t}$}
  \put(0.1,1.5){$\overline{s}$}
      \put(0.7,0.3){$\Delta_{\overline{t}, \overline{u}, \overline{s}}$}
  \end{picture}
}
\end{equation}
The slopes of edges of this triangle are $\overline{t}$, $\overline{u}$, $\overline{s}$ respectively.

\begin{Definition}[Ebina~\cite{Ebina}, Matsuura~\cite{Matsu}, Uchisawa~\cite{Uchi}]\label{Keu}
\begin{rm}
For $i = 1, 2, \ldots, u$, we put
\[
\ell_i = ^\# \{ (\alpha, \beta) \in \Delta_{\overline{t}, \overline{u}, \overline{s}} \cap {\Bbb Z}^2 \mid \alpha = i \} .
\]
Note that $\ell_u = 1$ and $\ell_i \ge 1$ for all $i = 1, 2, \ldots, u$.
We sort the sequence $\ell_1$, $\ell_2$, \ldots, $\ell_u$ into ascending order
\[
\ell'_1 \le \ell'_2 \le \cdots \le \ell'_u .
\]

We say that the {\em condition EMU} is satisfied for $(a,b,c)$ if
\[
\ell'_i \ge i
\]
for $i = 1, 2, \ldots, u$.
\end{rm}
\end{Definition}

Let me give an example.
Suppose $(a,b,c) = (17,503,169)$.
Then
\[
{\frak p} = (x^{89} - y^2z^3, y^{3}-x^{49}z^4,  z^7 - x^{40}y) 
\]
and the proper transform of $z^7 - x^{40}y$ is the negative curve.
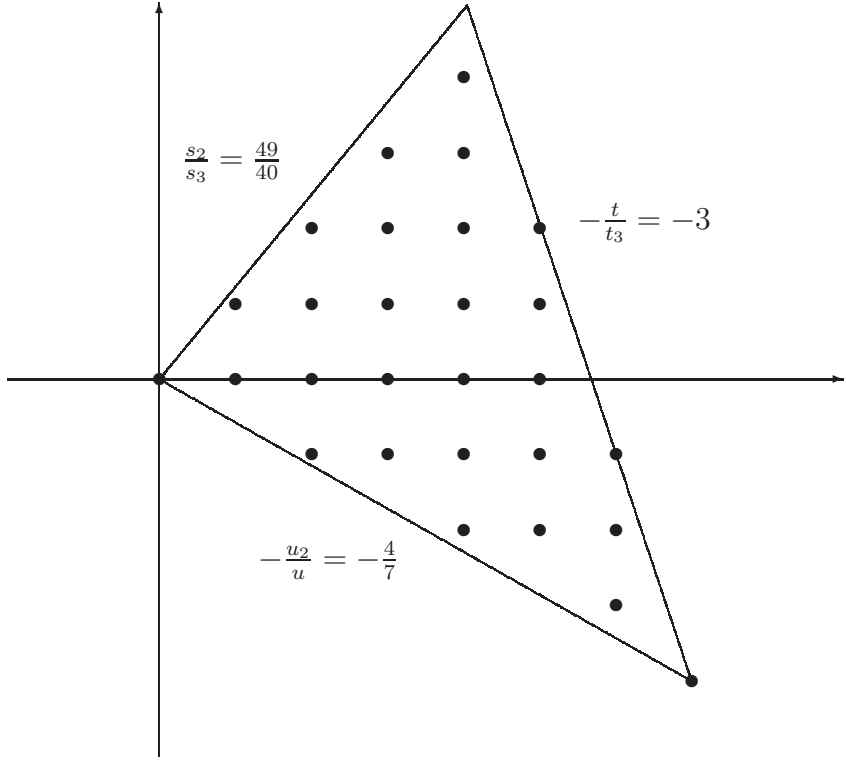
\[
{
\setlength\unitlength{1truecm}
  \begin{picture}(11,10)(-2,-5)
  \put(-2,0){\vector(1,0){11}}
  \put(0,-5){\vector(0,1){10}}
\qbezier (0,0) (3.5,-2) (7,-4)
\qbezier (0,0) (2.0238,2.4792) (4.0476,4.9583)
\qbezier (7,-4) (5.5238,0.47915)  (4.0476,4.9583)
\put(-0.1,-0.1){$\bullet$}
\put(0.9,-0.1){$\bullet$}
\put(0.9,0.9){$\bullet$}
\put(1.9,-1.1){$\bullet$}
\put(1.9,-0.1){$\bullet$}
\put(1.9,0.9){$\bullet$}
\put(1.9,1.9){$\bullet$}
\put(2.9,-1.1){$\bullet$}
\put(2.9,-0.1){$\bullet$}
\put(2.9,0.9){$\bullet$}
\put(2.9,1.9){$\bullet$}
\put(2.9,2.9){$\bullet$}
\put(3.9,-2.1){$\bullet$}
\put(3.9,-1.1){$\bullet$}
\put(3.9,-0.1){$\bullet$}
\put(3.9,0.9){$\bullet$}
\put(3.9,1.9){$\bullet$}
\put(3.9,2.9){$\bullet$}
\put(3.9,3.9){$\bullet$}
\put(4.9,-2.1){$\bullet$}
\put(4.9,-1.1){$\bullet$}
\put(4.9,-0.1){$\bullet$}
\put(4.9,0.9){$\bullet$}
\put(4.9,1.9){$\bullet$}
\put(5.9,-3.1){$\bullet$}
\put(5.9,-2.1){$\bullet$}
\put(5.9,-1.1){$\bullet$}
\put(6.9,-4.1){$\bullet$}
  \put(1.3,-2.5){$-\frac{u_2}{u} = -\frac{4}{7}$}
  \put(5.5,2){$-\frac{t}{t_3} = -3$}
  \put(0.3,2.8){$\frac{s_2}{s_3} = \frac{49}{40}$}
  \end{picture}
}
\]
Then $u = 7$ and 
\[
\ell_1 = 2, \ \ \ell_2 = 4, \ \ \ell_3 = 5, \ \ \ell_4 = 7, \ \ \ell_5 = 5, \ \ \ell_6 = 3, \ \ \ell_7 = 1 .
\]
Therefore
\[
\ell'_1 = 1, \ \ \ell'_2 = 2, \ \ \ell'_3 = 3, \ \ \ell'_4 = 4, \ \ \ell'_5 = 5, \ \ \ell'_6 = 5, \ \ \ell'_7 = 7 .
\]
The condition EMU is not satisfied in this case.

\vspace{2mm}

We put ${\frak p}^{(n)} = {\frak p}^nP_{\frak p}\cap P$, and call it the $n$th symbolic power of ${\frak p}$.
Consider the symbolic Rees ring
\[
R_s({\frak p}) := P[{\frak p}t, {\frak p}^{(2)}t^2, {\frak p}^{(3)}t^3, \ldots]
\subset P[t] .
\]

Assume that $z^u-x^{s_3}y^{t_3}$ is the negative curve, i.e., $\sqrt{abc} > uc$.
If the condition EMU is satisfied for $(a,b,c)$,
the symbolic Rees ring of ${\frak p}$ is Noetherian by Proposition~4.6 in \cite{KN}.
In this paper, we shall prove the converse as follows:

\begin{Theorem}\label{Keuprop}
Let $a$, $b$, $c$ be pairwise coprime positive integers.
Let $K$ be a field of characteristic $0$.
Assume that ${\frak p}$ is not complete intersection.
Suppose that $z^u-x^{s_3}y^{t_3}$ is a negative curve, i.e., $\sqrt{abc} > uc$.

Then $R_s({\frak p})$ is Noetherian if and only if 
the condition EMU is satisfied for $(a,b,c)$.
\end{Theorem}

It is an affirmative answer to a conjecture due to He~\cite{He}.
It was also implicitly conjectured in Kurano-Nishida~\cite{KN}.

\begin{Remark}
\begin{rm}
By the same argument as in the proof of Theorem~\ref{Keuprop},
we can prove the following:

Let $K$ be a field of characteristic $0$.
Let $\Delta_{r_1,r_2,r_3}$ be the triangle defined in Section~2.
We put $X_\Delta = {\rm Proj}(E(\Delta_{r_1, r_2, r_3}))$, where
$E(\Delta_{r_1, r_2, r_3})$ is the Ehrhart ring as in 
(\ref{EhrhartRing}).
Let $Y_\Delta$ be the blow-up of $X_\Delta$ at $(1,1) \in (K^\times)^2 \subset
X_\Delta$.
Assume
\[
\frac{1}{r_3-r_2}+\frac{1}{r_2-r_1} < 1 .
\]
Then the Cox ring of $Y_\Delta$ is Noetherian if and only if 
the condition EMU is satisfied for $\Delta_{r_1,r_2,r_3}$.
\end{rm}
\end{Remark}

The organization of this paper is as follows.

Section~2 is devoted to the preliminary.
We give a new criterion for finite generation in the case where there exists a negative curve $C$ such that $C.E = 1$ in section~3.
We give an algebraic description of this new criterion in section~4.
We shall give another proof to the finite generation in the case of ${\rm ch}(K) >0$ (Cutkosky~\cite{C}) and a result in the case of ${\rm ch}(K) =0$ (Kurano-Nishida~\cite{KN})
in section~5.
We classify $(a, b, c)$'s for which the condition EMU is not satisfied in section~6.
We prove Theorem~\ref{Keuprop} in section~7.

\section{Preliminary}\label{preliminary}

In this paper we assume that rings are commutative with $1$.
Let $\bN$ (resp.\ $\bNo$, $\bZ$, $\bR$, $\bRo$) be the set of positive integers (resp.\ non-negative integers, integers, real numbers, non-negative real numbers).

Let $r_1$, $r_2$, $r_3$ be rational numbers such that 
\begin{equation}\label{ri}
r_1 < -1 < r_2 < 0 < r_3 .
\end{equation}
Let $r_{ij}$ be integers such that 
\begin{equation}\label{rij}
\mbox{$r_i = r_{i1}/r_{i2}$ with
${\rm GCD}(r_{i1}, r_{i2})=1$ and $r_{i2}>0$ for $i = 1,2,3$.}
\end{equation}

Let $\Delta_{r_1,r_2,r_3}$ be the rational triangle with slopes $r_1$, $r_2$, $r_3$
such that its lower edge is the line segment connecting the end points $(0,0)$ and $(r_{22}, r_{12})$ as below.
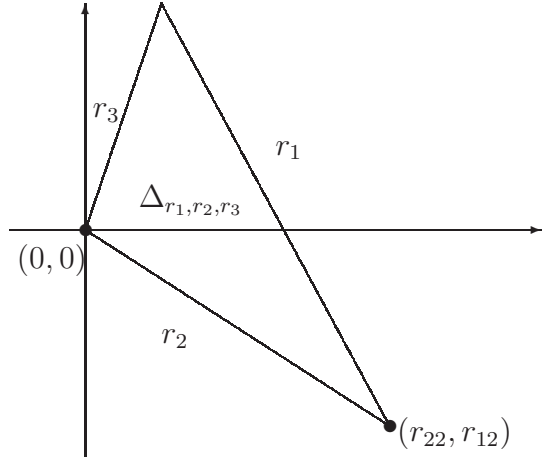
\[
{
\setlength\unitlength{1truecm}
  \begin{picture}(7,6)(-2,-3)
  \put(-1,0){\vector(1,0){7}}
  \put(0,-3){\vector(0,1){6}}
\qbezier (0,0) (2,-1.3) (4,-2.6)
\qbezier (0,0) (0.5,1.5) (1,3)
\qbezier (1,3) (2.5,0.2) (4,-2.6)
\put(-0.1,-0.1){$\bullet$}
\put(3.9,-2.7){$\bullet$}
\put(-0.9,-0.5){$(0,0)$}
\put(4.1,-2.8){$(r_{22}, r_{12})$}
  \put(1,-1.5){$r_2$}
  \put(2.5,1){$r_1$}
  \put(0.1,1.5){$r_3$}
      \put(0.7,0.3){$\Delta_{r_1,r_2,r_3}$}
  \end{picture}
}
\]
We define the Ehrhart ring of $\Delta_{r_1, r_2, r_3}$ over a field $K$ by 
\begin{equation}\label{EhrhartRing}
E(\Delta_{r_1, r_2, r_3})
= \bigoplus_{(\alpha,\beta) \in m\Delta_{r_1, r_2, r_3} \cap \bZ^2}
Kv^{\alpha}w^\beta t^m \subset K[v^{\pm 1}, w^{\pm 1}, t] ,
\end{equation}
where $v$, $w$, $t$ are algebraically independent over $K$.

Let $K$ be a field.
Suppose that $a$, $b$, $c$ are pairwise coprime positive integers.
Let ${\frak p}$ be the kernel of the $K$-algebra map
$\varphi : P = K[x,y,z] \rightarrow K[T]$ defined by 
$\varphi(x) = T^a$, $\varphi(y) = T^b$, $\varphi(z) = T^c$.
In this paper we always assume that 
\[
\mbox{${\frak p}$ is not complete intersection.}
\]
Then we know
\[
{\frak p} = I_2\left(
\begin{array}{ccc}
x^{s_2} & y^{t_3} & z^{u_1} \\
y^{t_1} & z^{u_2} & x^{s_3}
\end{array}
\right) = (x^s-y^{t_1}z^{u_1}, y^t-x^{s_2}z^{u_2}, z^u-x^{s_3}y^{t_3})
\]
with positive integers $s_2$, $s_3$, $t_1$, $t_3$, $u_1$, $u_2$,
where $I_2( \ )$ is the ideal generated by $2 \times 2$ minors of this given $2 \times 3$-matrix, and
 $s = s_2+s_3$, $t = t_1+t_3$, $u = u_1+u_2$ (Herzog~\cite{Her}).
We put, 
$\overline{t} = -t/t_3$, $\overline{u}=-u_2/u$, $\overline{s}=s_2/s_3$.
We think that $P=K[x,y,z]$ is a graded polynomial ring with ${\rm deg}(x) = a$, ${\rm deg}(y) = b$, ${\rm deg}(z) = c$.
Then the Veronesean subring $S^{(ab)}$ is isomorphic to $E(\Delta_{\overline{t}, \overline{u}, \overline{s}})$ as in section~4 in \cite{KN}.

We have the following proposition:

\begin{Proposition}\label{torsion}
Let $r_i$ and $r_{ij}$ be numbers satisfying (\ref{ri}) and (\ref{rij}).
Then the following three conditions are equivalent:
\begin{enumerate}
\item
The class group of $\proj E(\Delta_{r_1,r_2,r_3})$ is torion-free.
\item
$\bZ \left( \begin{array}{c} r_{11} \\ r_{12} \end{array} \right) + \bZ \left( \begin{array}{c} r_{21} \\ r_{22}  \end{array} \right) + \bZ \left( \begin{array}{c} r_{31} \\ r_{32}  \end{array} \right) =\bZ^2$. 
\item
There exist pairwise coprime positive integers $a$, $b$, $c$ such that
$\overline{t} = r_1$, $\overline{u} = r_2$ and $\overline{s} = r_3$.
\end{enumerate}
\end{Proposition}

We omit a proof.
By this proposition, we immediately obtain the following corollary.

\begin{Corollary}\label{abc}
Let $r_i$ and $r_{ij}$ be numbers satisfying (\ref{ri}) and (\ref{rij}).
Let $\epsilon$ be any positive number.
Then there exists a positive number $r'_3$ with $r_3-\epsilon < r'_3 < r_3$ such that
there exist pairwise coprime positive integers $a$, $b$, $c$ such that
$\overline{t} = r_1$, $\overline{u} = r_2$ and $\overline{s} = r'_3$.
\end{Corollary}

\begin{Remark}
\begin{rm}
There exists a one to one correspondence between
\[
A:= \left\{ (a,b,c) \ \left| \ 
\begin{tabular}{l}
\mbox{$a$, $b$, $c$ are pairwise coprime positive integers,} \\ 
\mbox{${\frak p}$ is not complete intersection}
\end{tabular}
\right\}
\right.
\]
and
\[
B:= \left\{ \Delta_{r_1,r_2,r_3} \ \left| \ 
\begin{tabular}{l}
\mbox{$r_1$, $r_2$, $r_3$ are rational numbers,} \\
\mbox{$r_1 < -1 < r_2 < 0 < r_3$,}  \\ 
\mbox{$\bZ \left( \begin{array}{c} r_{11} \\ r_{12} \end{array} \right) + \bZ \left( \begin{array}{c} r_{21} \\ r_{22}  \end{array} \right) + \bZ \left( \begin{array}{c} r_{31} \\ r_{32}  \end{array} \right) =\bZ^2$}
\end{tabular}
\right\}
\right. .
\]

We define a map $\Psi:A \rightarrow B$ to be $\Psi((a,b,c)) = \Delta_{\overline{t}, \overline{u}, \overline{s}}$.

We define a map $\Phi:B \rightarrow A$ to be 
$\Phi(\Delta_{r_1,r_2,r_3}) = (a,b,c)$, where $a$, $b$, $c$ are pairwise coprime positive integers satisfying 
\[
b \left( \begin{array}{c} r_{11} \\ r_{12} \end{array} \right) - c \left( \begin{array}{c} r_{21} \\ r_{22}  \end{array} \right) +a \left( \begin{array}{c} r_{31} \\ r_{32}  \end{array} \right) =
\left( \begin{array}{c} 0 \\ 0  \end{array} \right) .
\]
By Proposition~\ref{torsion}, we can prove that $\Phi$ is the inverse mapping of $\Psi$.
\end{rm}
\end{Remark}

%%%%%%%%%%%%%%%%%%%%%%%%%%%%%%%%%%%%%%%%%%%%%%%%%%%%%%%%%%%%%%%%%%%
\section{A new criterion for finite generation of symbolic Rees rings}

Let $K$, $a$, $b$, $c$, $P=K[x,y,z]$, $\bfp$, $\Delta_{\overline{t}, \overline{u}, \overline{s}}$ be as in section~\ref{preliminary}.
In the rest of this paper, we always assume that
\begin{equation}\label{abcjyouken}
\begin{split}
\bullet \ \ & 
\mbox{${\frak p}$ is not complete intersection, and} \\
\bullet \ \ & 
\mbox{$z^u-x^{s_3}y^{t_3}$ is the negative curve, i.e., $uc < \sqrt{abc}$.}
\end{split}
\end{equation}
The second condition above is equivalent to that
the area of $\Delta_{\overline{t}, \overline{u}, \overline{s}}$ is bigger than $u^2/2$.
It is also equivalent to 
\begin{equation}\label{stu}
\frac{1}{\overline{s}-\overline{u}} + \frac{1}{\overline{u}-\overline{t}}<1 .
\end{equation}

\[
{
\setlength\unitlength{1truecm}
  \begin{picture}(7,6)(0,-3)
  %\put(-1,0){\vector(1,0){7}}
  %\put(0,-3){\vector(0,1){6}}
\qbezier (0,0) (3,-1.5) (6,-3)
\qbezier (3,-1.5) (2,0.5) (1,2.5)
\qbezier (3,-1.5) (4,1) (5, 3.5)
\put(2.9,-1.6){$\bullet$}
\put(2.5,-2){$(0,0)$}
\put(4,-1){$S$}
  \put(1,-1){$\overline{u}$}
  \put(1.9,1){$\overline{t}$}
  \put(3.8,1.5){$\overline{s}$}
      \put(0.7,0.3){$T$}
  \end{picture}
}
\]
%
%\[
%\begin{tikzpicture}
%\tikz \filldraw[fill=lightgray,very thick] (3,-1)--(0,0)--(2,4);
% \draw (-0.5,1)node{$(0,0)$};
% \draw (1.5,0.2)node{$\overline{u}$};
% \draw (1,3.5)node{$\overline{s}$};
%  \draw (1,1.5)node{$S$};
%\end{tikzpicture}
%\ \ \ \ \ 
%\begin{tikzpicture}
%\tikz \filldraw[fill=lightgray,very thick] (-3,1)--(0,0)--(-3,5);
% \draw (3.7,0)node{$(0,0)$};
% \draw (1.5,0.2)node{$\overline{u}$};
% \draw (1.5,3.5)node{$\overline{t}$};
%  \draw (1,1.5)node{$T$};
%\end{tikzpicture}
%\]
Let $S$ and $T$ be the cone in $\bR^2$ defined by 
\begin{align*}
S & = \bRo (u,-u_2) + \bRo (s_3,s_2),  \\
T & = \bRo (-u,u_2) + \bRo (-t_3,t) .
\end{align*}
We put
\begin{align*}
K[S] & = \bigoplus_{(\alpha,\beta) \in S \cap \bZ^2} Kv^\alpha w^\beta \subset K[v^{\pm 1}, w^{\pm 1}] , \\
K[T] & = \bigoplus_{(\alpha,\beta) \in T \cap \bZ^2} Kv^\alpha w^\beta \subset K[v^{\pm 1}, w^{\pm 1}] .
\end{align*}

Let $\pi : Y \rightarrow X = \proj P$ be the blow-up at $V_+(\bfp)$.
Let $C$ be the proper transform of $V_+(z^u-x^{s_3}y^{t_3})$.
Then we have $\pi(C) = V_+(z^u-x^{s_3}y^{t_3}) \simeq \bP^1_K$.
Put $U_1 = \spec K[S]$ and $U_2 = \spec K[T]$.
Then $U_i$'s are affine open sets of $X$ such that $U_1\cup U_2 \supset \pi(C)$.
Since $C$ is isomorphic to $\pi(C)$, 
\[
C = (\pi^{-1}(U_1) \cap C) \cup (\pi^{-1}(U_2) \cap C)
\]
is an affine open covering of $C$.
For a positive integer $\ell$, $\ell C$ is the closed subscheme of $Y$ defined by
the ideal sheaf $\oo_Y(-\ell C)$.
In general, for a scheme $W$, $W$ is affine if and only if so is $W_{\rm red}$ 
(e.g., Exercise~3.1 in section~\RMN{3} in \cite{H}).
Therefore
\begin{equation}\label{nCaffine}
\ell C = (\pi^{-1}(U_1) \cap \ell C) \cup (\pi^{-1}(U_2) \cap \ell C)
\end{equation}
is an affine open covering of $\ell C$.

Let $H$ be a Weil divisor on $Y$ satisfying $\oo_Y(H) = \pi^*\oo_X(1)$.
(Here remark that $V_+(\bfp)$ is a non-singular point of $X$.
Therefore $\pi^*\oo_X(1)$ is a reflexive sheaf on $Y$.)
Let $E$ be the exceptional divisor of the blow-up $\pi : Y \rightarrow X$.
Then we have 
\[
{\rm Cl}(Y) = \bZ H + \bZ E \simeq \bZ^2
\]
and 
\[
H^2 = \frac{1}{abc}, \ \ E^2 = -1, \ \ H.E = 0.
\]

Then we have the following criterion for finite generation of $R_s({\frak p})$.

\begin{Theorem}[Cutkosky]\label{Cut}
The following conditions are equivalent:
\begin{enumerate}
\item
$R_s({\frak p})$ is finitely generated over $P$ (equivalently,  $R_s({\frak p})$ is Noetherian).
\item
The Cox ring ${\rm Cox}(Y)$ of $Y$ is finitely generated over $K$ (equivalently,  ${\rm Cox}(Y)$ is Noetherian).
\item
There exists a curve $D$ on $Y$ satisfying $C\cap D = \emptyset$.
\end{enumerate}
\end{Theorem}

Since $C \sim ucH-E$, we know
$C.(abH-uE) = 0$.
One can show that, if there exists a curve $D$ satisfying the condition (3) in Theorem~\ref{Cut}, there exists a positive integer $m$ such that 
$D \sim m(abH-uE)$.
(Since $D$ does not contain the points corresponding to the end points of the lower edge of $\Delta_{\overline{t}, \overline{u}, \overline{s}}$, 
we know that the degree of $\pi(D)$ must be a multiple of $ab$.)
The triangle $\Delta_{\overline{t}, \overline{u}, \overline{s}}$ is corresponding to the toric divisor $\oo_X(ab)$.
We know that $abH-uE$ is a Cartier divisor at each point in $Y$ except for the point corresponding to the top vertex of $\Delta_{\overline{t}, \overline{u}, \overline{s}}$.
On the other hand, the curve $C$ does not pass through this point.
Ultimately, $\oo_Y(m(abH-uE))|_{\ell C}$ is an invertible sheaf on $\ell C$ for any $m$ and $\ell$.

Here we have the following theorem:

\begin{Theorem}\label{new}
The following conditions are equivalent:
\begin{enumerate}
\item
$R_s({\frak p})$ is finitely generated over $P$.
\item
There exists a positive integer $m$ such that $\oo_Y(m(abH-uE))|_{\ell C} \simeq \oo_{\ell C}$ for any positive integer $\ell$.
\item
There exists a positive integer $m$ such that $\oo_Y(m(abH-uE))|_{muC} \simeq \oo_{muC}$.
\end{enumerate}
\end{Theorem}

\proof
First we shall show $(1) \Longrightarrow (2)$.
Since $R_s({\frak p})$ is finitely generated over $P$, there exists a curve $D$ on $Y$ satisfying $C\cap D = \emptyset$ as in Theorem~\ref{Cut}.
Then $D \sim m(abH-uE)$ for some positive integer $m$.
Since $C \cap D = \emptyset$, the section in $H^0(Y, \oo_Y(m(abH-uE)))$ corresponding to $D$ does not vanish at any point in $C$.
Hence it does not vanish at any point in $\ell C$ for any positive integer $\ell$.
Therefore $\oo_Y(m(abH-uE))|_{\ell C}$ has a global sention nowhere vanishes.
Thus we have obtained $\oo_Y(m(abH-uE))|_{\ell C} \simeq \oo_{\ell C}$.

The implication $(2) \Longrightarrow (3)$ is trivial.

Next we shall prove $(3) \Longrightarrow (1)$.
Consider the tensor product
\[
\oo_Y(m(abH-uE)) \otimes_{\oo_Y}
\left( 0 \rightarrow \oo_Y(-muC) \rightarrow \oo_Y \rightarrow \oo_{muC} \rightarrow 0 \right) .
\]
Here remark that $\oo_Y(m(abH-uE))$ is invertible at any point in $muC$.
Therefore we have the following exact sequence
\[
0 \rightarrow \oo_Y(m(abH-uE-uC)) \rightarrow \oo_Y(m(abH-uE)) \rightarrow \oo_Y(m(abH-uE))|_{muC} \rightarrow 0 .
\]
Taking a long exact sequence on cohomologies, we obtain
the following exact sequence:
\begin{equation}\label{mapf}
H^0(Y,\oo_Y(m(abH-uE))) \stackrel{\psi}{\rightarrow} H^0(mu C, \oo_Y(m(abH-uE))|_{muC})
\rightarrow H^1(Y,\oo_Y(m(abH-uE-uC))) .
\end{equation}
Since $C \sim ucH-E$, we know
\[
H^1(Y,\oo_Y(m(abH-uE-uC))) = H^1(Y,\oo_Y(m(ab-u^2c)H)) .
\]
Since 
\[
R^i\pi_*(\oo_Y(m(ab-u^2c)H)) 
= \left\{
\begin{array}{ll}
\oo_X(m(ab-u^2c)) & (i=0) \\
0 & (i>0)
\end{array}
\right. ,
\]
we have
\[
H^1(Y,\oo_Y(m(ab-u^2c)H)) \simeq H^1(X,\oo_X(m(ab-u^2c))) 
\]
by the Leray spectral sequence.
Then we have
\[
H^1(X,\oo_X(m(ab-u^2c))) \simeq H^2_{(x,y,z)}(P)_{m(ab-u^2c)} = 0 .
\]
Ultimately, the map $\psi$ in the exact sequence (\ref{mapf}) is surjective.
Since $\oo_Y(m(abH-uE))|_{muC} \simeq \oo_{muC}$,
there exists $f \in H^0(Y,\oo_Y(m(abH-uE)))$ such that
$\psi(f)$ does not vanish at any point in $muC$.
Then any irreducible component corresponding to the effective divisor ${\rm div}(f) + m(abH-uE)$ does not meet $C$.
Therefore $R_s({\frak p})$ is Noetherian by Theorem~\ref{Cut}.
\qed

\begin{Corollary}\label{mwofuyasu}
Assume $\oo_Y(m(abH-uE))|_{muC} \simeq \oo_{muC}$ for some $m>0$.
Then we have $\oo_Y(mm'(abH-uE))|_{mm'uC} \simeq \oo_{mm'uC}$ for any $m'>0$.
\end{Corollary}

\proof
By Theorem~\ref{new}, we have $\oo_Y(m(abH-uE))|_{mm'uC} \simeq \oo_{mm'uC}$ for any $m'>0$.
The assertion immediately follows from it.
\qed

\vspace{2mm}

Applying this theorem,
we can give another proof to 
\begin{itemize}
\item
$\oo_Y(p^e(abH-uE))|_{p^euC} \simeq \oo_{p^euC}$ for $e\gg 0$, in particular $R_s({\frak p})$ is always Noetherian in the case of ${\rm ch}(K)=p > 0$ (Cutkosky~\cite{C}), and
\item
that the finite generation of $R_s({\frak p})$ is equivalent to 
$\oo_Y(abH-uE)|_{uC} \simeq \oo_{uC}$ in the case of ${\rm ch}(K) = 0$ (Kurano-Nishida~\cite{KN})
\end{itemize}
in Proposition~\ref{ch=p>0} and \ref{ch=0}.

\section{Algebraic description of Theorem~\ref{new}}

In this section we shall describe Theorem~\ref{new} using algebraic method.

Consider the affine openset $U_1 = \spec K[S]$ of $X$.
Then
\[
M_1 = (v-1) + (\{ v^\alpha w^\beta (w-1) \mid (\alpha,\beta), (\alpha,\beta+1) \in S \} )
\]
is the maximal ideal of $K[S]$ corresponding to $V_+({\frak p})$.
Consider the following affine openset of the blow-up at $M_1$:
\[
A := K[S]\left[ \left\{ \left. \frac{v^\alpha w^\beta (w-1)}{v-1} \ \right| \
(\alpha,\beta), (\alpha,\beta+1) \in S \right\} \right]
\subset K[v^{\pm 1}, w^{\pm 1}, \frac{w-1}{v-1}] .
\]
Then $\spec A$ is an affine open subset of $Y$.
The defining ideal of $C$ is
\[
\left( \left\{ \left. \frac{v^\alpha w^\beta (w-1)}{v-1} \ \right| \
(\alpha,\beta), (\alpha,\beta+1) \in S \right\} \right)A .
\]

Consider the affine openset $U_2 = \spec K[T]$ of $X$.
Then
\[
M_2 = (v^{-1}w-1) + (\{ v^\alpha w^\beta (w-1) \mid (\alpha,\beta), (\alpha,\beta+1) \in T \} )
\]
is the maximal ideal of $K[T]$ corresponding to $V_+({\frak p})$.
Here remark $\overline{t} < -1 < \overline{u} < 0$ and $v^{-1}w-1 \in K[T]$.
Consider the following affine openset of the blow-up at $M_2$:
\[
B := K[T]\left[ \left\{ \left. \frac{v^\alpha w^\beta (w-1)}{v^{-1}w-1} \ \right| \
(\alpha,\beta), (\alpha,\beta+1) \in T \right\} \right]
\subset K[v^{\pm 1}, w^{\pm 1}, \frac{w-1}{v^{-1}w-1}] .
\]
Then $\spec B$ is an affine open subset of $Y$.
The defining ideal of $C$ is
\[
\left( \left\{ \left. \frac{v^\alpha w^\beta (w-1)}{v^{-1}w-1} \ \right| \
(\alpha,\beta), (\alpha,\beta+1) \in T \right\} \right)B .
\]

Here it is easy to see 
\[
\pi^{-1}(U_1) \cap \ell C \subset \spec A, \ \ 
\pi^{-1}(U_2) \cap \ell C \subset \spec B 
\]
for any $\ell>0$.
By (\ref{nCaffine}),
\begin{equation}\label{coverofellC}
\ell C = (\spec A)|_{\ell C} \cup (\spec B)|_{\ell C}
\end{equation}
is an affine open covering of $\ell C$.

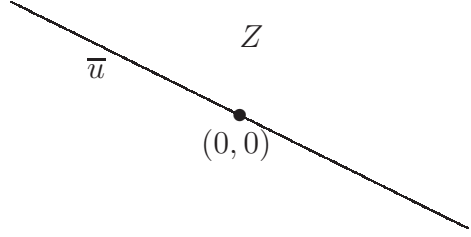
\[
{
\setlength\unitlength{1truecm}
  \begin{picture}(7,4)(0,-3)
  %\put(-1,0){\vector(1,0){7}}
  %\put(0,-3){\vector(0,1){6}}
\qbezier (0,0) (3,-1.5) (6,-3)
%\qbezier (3,-1.5) (2,0.5) (1,2.5)
%\qbezier (3,-1.5) (4,1) (5, 3.5)
\put(2.9,-1.6){$\bullet$}
\put(2.5,-2){$(0,0)$}
%\put(4,-1){$S$}
  \put(1,-1){$\overline{u}$}
  %\put(2,1){$\overline{t}$}
  %\put(3.7,1.5){$\overline{s}$}
      %\put(0.7,1){$T$}
      \put(3,-0.6){$Z$}
  \end{picture}
}
\]
%\[
%\begin{tikzpicture}
%\draw[line width=1.5pt] (-2,0.5)--(2,-0.5);
%\coordinate (O) at (0,0);
% \foreach \P in {O} \fill[black] (\P) circle (0.06); 
%  \coordinate  (A) at (-2,0.5); %点A
% \coordinate  (B) at (2,-0.5); %点B
% \coordinate  (C) at (2,2);
%  \coordinate  (D) at (-2,2);
%  \fill[lightgray] (A)--(B)--(C)--(D)--cycle;
%   \draw (-1,0)node{$\overline{u}$};
%  \draw (0,1)node{$Z$};
%    \draw (0,-0.5)node{$(0,0)$};
%\end{tikzpicture}
%\]
Let $Z$ be the cone $\bR(u,-u_2) + \bRo(0,1)$ as above.
Put
\[
x = \frac{w-1}{v-1}, \ \ y = \frac{w-1}{v^{-1}w-1} .
\]
Put 
\begin{align*}
K[Z] & = \bigoplus_{(\alpha,\beta) \in Z \cap \bZ^2} Kv^\alpha w^\beta \subset K[v^{\pm 1}, w^{\pm 1}] , \\
L & = K[Z]\left[ x, y,
\frac{v-1}{v^{-1}w-1}, \frac{v^{-1}w-1}{v-1} \right] 
\subset K[v^{\pm 1},w^{\pm 1}]\left[ \frac{1}{v-1}, \frac{1}{v^{-1}w-1} \right] , \\
F & = K[Z]\left[ x \right] \subset L , \\
G & = K[Z]\left[ y \right] \subset L .
\end{align*}
Here remark that
$v$, $w$, $v^{-1}w$ are contained in $K[Z]$ since $-1 < \overline{u} < 0$.
Since 
\begin{equation}\label{unit}
\frac{v^{-1}w-1}{v-1} = -v^{-1}w + x \in F ,
\end{equation}
the map $F \hookrightarrow L$ is a localization.
Since 
\[
\frac{v-1}{v^{-1}w-1} = -v + y \in G ,
\]
the map $G \hookrightarrow L$ is a localization.

Let $\ell$ be a positive integer.
We know that $w$ is a unit of $F/x^\ell F$ since $F/xF = K[Z]/(w-1)$.
We know that $v^{-1}w$ is a unit of $F/x^\ell F$ since $w = (v^{-1}w)v$.
Since $x$ is a nilpotent in $F/x^\ell F$, 
$\frac{v^{-1}w-1}{v-1}$ is a unit in $F/x^\ell F$ by (\ref{unit}).
We can show that $\frac{v-1}{v^{-1}w-1}$ is a unit in $G/y^\ell G$ in the same way.
Thus we obtain isomorphisms
\begin{equation}\label{FtoG}
F/x^\ell F \simeq L/x^\ell L = L/y^\ell L \simeq G/y^\ell G
\end{equation}
for any $\ell >0$.

\vspace{2mm}

For $\alpha \in \bZ$ and $n \ge 0$, we define
\[
x_{\alpha,n} := v^\alpha w^{\lceil \alpha \overline{u} \rceil} x^n \in F ,
\]
where $\lceil \alpha \overline{u} \rceil$ is the least integer such that $\lceil \alpha \overline{u} \rceil \ge \alpha \overline{u}$.

\begin{Proposition}\label{strF}
We have 
\begin{align}
\label{strF1} F & = \bigoplus_{\alpha \in \bZ} \bigoplus_{n \ge 0} Kx_{\alpha,n} \\
\label{strF2}  x^\ell F & = \bigoplus_{\alpha \in \bZ} \bigoplus_{n \ge \ell} Kx_{\alpha,n} 
\end{align}
for any positive integer $\ell$.
\end{Proposition}

\proof
We have
\[
F = K[Z][x] = 
\sum_{(\alpha,\beta) \in Z\cap \bZ^2} \sum_{n\ge 0}K v^\alpha w^\beta x^n
\supset 
\sum_{\alpha \in \bZ} \sum_{n\ge 0}K x_{\alpha,n} 
\]
by definition.

We shall prove the opposite containment.
For rational numbers $p$, $q$, $r$ such that $q<r$, let $W_{p,q,r}$ be the triangle
with vertices ${\rm A}(p,q)$, ${\rm B}(p,r)$, ${\rm C}(\frac{r+p-q+p\overline{u}}{1+\overline{u}},\frac{r+q\overline{u}}{1+\overline{u}})$. 
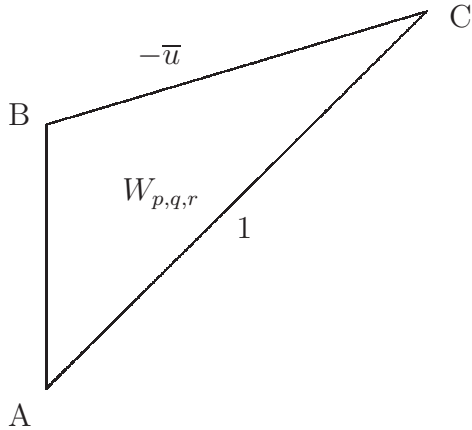
\[
{
\setlength\unitlength{1truecm}
  \begin{picture}(7,6)(1,-2)
  %\put(-1,0){\vector(1,0){7}}
  %\put(0,-3){\vector(0,1){6}}
\qbezier (3,2) (5.5,2.75) (8,3.5)
\qbezier (3,-1.5) (3,1) (3,2)
\qbezier (3,-1.5) (4,-0.5) (8,3.5)
%\put(2.9,-1.6){$\bullet$}
\put(2.5,-2){A}
\put(2.5,2){B}
\put(8.3,3.3){C}
  \put(5.5,0.5){$1$}
  \put(4.2,2.8){$-\overline{u}$}
  \put(4,1){$W_{p,q,r}$}
      %\put(0.7,0.3){$T$}
  \end{picture}
}
\]
%\[
%\begin{tikzpicture}
% \coordinate[label=left:A] (A) at (0,0); 
% \coordinate[label=right:C] (C) at (4,4); 
% \coordinate[label=left:B] (B) at (0,2.5); 
% \fill[lightgray] (A)--(B)--(C)--cycle; 
% \draw[very thick] (A)--(B)--(C)--cycle; 
%    \draw (1.8,3.5)node{$-\overline{u}$};
%        \draw (2.5,2.2)node{$1$};
%  \draw (1,1.8)node{$W_{p,q,r}$};
%\end{tikzpicture}
%\]
The slope of the lower edge of $W_{p,q,r}$ is $1$, and that of the upper edge is $-\overline{u}$.
We shall prove the following claim:

\begin{Claim}\label{claim}
Let $\alpha$, $\beta$, $n$ be integers such that $\beta \ge \alpha \overline{u}$ and $n\ge 0$.
Then we have\footnote{
If the characteristic of $K$ is $0$, we can prove
$v^\alpha w^\beta x^n \in \sum_{(\gamma, m) \in W_{\alpha,n,\beta-\alpha\overline{u}+n} \cap \bZ^2}\bZ x_{\gamma, m}$.
} 
\[
v^\alpha w^\beta x^n \in \sum_{(\gamma, m) \in W_{\alpha,n,\beta-\alpha\overline{u}+n} \cap \bZ^2}Kx_{\gamma, m} ,
\]
where the coefficient of $x_{\alpha,n}$ in $v^\alpha w^\beta x^n$ is $1$.
\end{Claim}

First suppose $0 \le \beta - \alpha\overline{u} < 1$.
Then $v^\alpha w^\beta x^n$ is equal to $x_{\alpha, n}$ since 
$\beta = \lceil \alpha \overline{u} \rceil$.

Next suppose $\beta - \alpha\overline{u} \ge 1$.
Since $w^\beta = w^{\beta-1}+w^{\beta-1}(w-1)$, we have
\begin{equation}\label{vwx}
v^\alpha w^\beta x^n = v^\alpha w^{\beta-1} x^n+v^\alpha w^{\beta-1}(v-1) x^{n+1}
=v^\alpha w^{\beta-1} x^n+v^{\alpha+1} w^{\beta-1} x^{n+1}-v^\alpha w^{\beta-1} x^{n+1} .
\end{equation}
We assume that this claim is true if $\beta - \alpha\overline{u}$ is smaller.
Then we know
\begin{align*}
v^\alpha w^{\beta-1} x^n & \in \sum_{(\gamma, m) \in W_{\alpha,n,\beta-1-\alpha\overline{u}+n} \cap \bZ^2}Kx_{\gamma, m} , \\
v^{\alpha+1} w^{\beta-1} x^{n+1} & \in \sum_{(\gamma, m) \in W_{\alpha+1,n+1,\beta-1-(\alpha+1)\overline{u}+n+1} \cap \bZ^2}Kx_{\gamma, m} , \\
v^\alpha w^{\beta-1} x^{n+1} & \in \sum_{(\gamma, m) \in W_{\alpha,n+1,\beta-1-\alpha\overline{u}+n+1} \cap \bZ^2}Kx_{\gamma, m} .
\end{align*}
It is easy to see that triangles $W_{\alpha,n,\beta-1-\alpha\overline{u}+n}$,
$W_{\alpha+1,n+1,\beta-1-(\alpha+1)\overline{u}+n+1}$ and
$W_{\alpha,n+1,\beta-1-\alpha\overline{u}+n+1}$ are contained in 
$W_{\alpha,n,\beta-\alpha\overline{u}+n}$.
Thus we have proved the claim by (\ref{vwx}).

In the rest of this proof, we shall show that $x_{\alpha,n}$'s are linearly independent over $K$.
Suppose
\begin{equation}\label{linindep}
\sum_{\alpha,n}c_{\alpha,n}x_{\alpha,n} = 0 \ \ (c_{\alpha,n}\in K) .
\end{equation}
Let $\gamma$ be the minimal number which satisfies $c_{\gamma,m} \neq 0$ for some $m\ge 0$.
We think
\[
x_{\alpha,n} \in F \subset K[[v,w]][v^{-1},w^{-1}] .
\]
Taking the coefficient of $v^\gamma$ of (\ref{linindep}),
we have
\[
\sum_{n\ge 0}c_{\gamma,n}w^{\lceil \gamma \overline{u} \rceil} (1-w)^n = 0 .
\]
It implies $c_{\gamma,n} = 0$ for any $n \ge 0$.
It is a contradiction.
We have proved (\ref{strF1}) in Proposition~\ref{strF}.

We have
\[
x^\ell F = 
\sum_{(\alpha,\beta) \in Z\cap \bZ^2} \sum_{n\ge \ell}K v^\alpha w^\beta x^n
\supset 
\sum_{\alpha \in \bZ} \sum_{n\ge \ell}K x_{\alpha,n} 
= \bigoplus_{\alpha \in \bZ} \bigoplus_{n\ge \ell}K x_{\alpha,n} .
\]
Using Claim~\ref{claim}, we can prove the opposite containment.
\qed

\begin{Remark}\label{algstr}
\begin{rm}
For integers $\alpha$, $n$, let $V_{\alpha,n}$ be the cone with vertex ${\rm A}(\alpha,n)$ as below:
\[
{
\setlength\unitlength{1truecm}
  \begin{picture}(7,6)(1,-2)
  %\put(-1,0){\vector(1,0){7}}
  %\put(0,-3){\vector(0,1){6}}
%\qbezier (3,2) (5.5,2.75) (8,3.5)
\qbezier (3,-1.5) (3,1) (3,3.5)
\qbezier (3,-1.5) (4,-0.5) (8,3.5)
%\put(2.9,-1.6){$\bullet$}
\put(2.5,-2){A}
%\put(2.5,2){B}
%\put(8.3,3.3){C}
  \put(5.5,0.5){$1$}
 % \put(4.2,2.8){$-\overline{u}$}
  \put(4,1){$V_{\alpha,n}$}
      %\put(0.7,0.3){$T$}
  \end{picture}
}
\]
%\[
%\begin{tikzpicture}
%\tikz \filldraw[fill=lightgray,very thick] (0,4)--(0,0)--(4,4);
% \coordinate[label=left:A] (A) at (0,0); %点Aの定義
%         \draw (2.5,2.2)node{$1$};
%  \draw (1,1.8)node{$V_{\alpha,n}$};
%\end{tikzpicture}
%\]

By Claim~\ref{claim}, we have\footnote{
If the characteristic of $K$ is $0$, we can prove
$x_{\alpha,n}x_{\alpha',n'} = v^{\alpha+\alpha'}w^{\lceil \alpha \overline{u} \rceil + \lceil \alpha' \overline{u} \rceil}x^{n+n'}\in \sum_{(\gamma,m) \in V_{\alpha+\alpha',n+n'} \cap \bZ^2}
\bZ x_{\gamma,m}$.
} 
\begin{equation}\label{algstr}
x_{\alpha,n}x_{\alpha',n'} = v^{\alpha+\alpha'}w^{\lceil \alpha \overline{u} \rceil + \lceil \alpha' \overline{u} \rceil}x^{n+n'}\in \sum_{(\gamma,m) \in V_{\alpha+\alpha',n+n'} \cap \bZ^2}
Kx_{\gamma,m} ,
\end{equation}
where the coefficient of $x_{\alpha+\alpha',n+n'}$ is $1$.
Here remark $\lceil (\alpha +  \alpha') \overline{u} \rceil \le \lceil \alpha \overline{u} \rceil + \lceil \alpha' \overline{u} \rceil$.
\end{rm}
\end{Remark}

\begin{Proposition}\label{strA}
We have the following equalities:
\begin{align}
\label{strA1}
A & = \sum_{\tiny \begin{array}{c} n \ge 0 \\ (\alpha,\beta), (\alpha,\beta+n) \in S \end{array}} Kv^\alpha w^\beta x^n = 
\bigoplus_{\alpha \ge 0} \bigoplus_{\tiny \begin{array}{c} n \ge 0 \\ (\alpha,\lceil \alpha \overline{u} \rceil+n) \in S \end{array}} Kx_{\alpha,n} \\
\label{strA2}
x^\ell F \cap A & = \sum_{\tiny \begin{array}{c} n \ge \ell  \\ (\alpha,\beta), (\alpha,\beta+n) \in S \end{array}} Kv^\alpha w^\beta x^n= \bigoplus_{\alpha \ge 0} \bigoplus_{\tiny \begin{array}{c} n \ge \ell  \\ (\alpha,\lceil \alpha \overline{u} \rceil+n) \in S \end{array}} Kx_{\alpha,n} 
\end{align}
for any positive integer $\ell$.
\end{Proposition}

\proof
First we shall prove
\[
A = \sum_{\tiny \begin{array}{c} n \ge 0 \\ (\alpha,\beta), (\alpha,\beta+n) \in S \end{array}} Kv^\alpha w^\beta x^n .
\]
It is easy to see that the right hand side is a subring of $K[v^{\pm 1}, w^{\pm 1}, \frac{1}{v-1}]$.
Therefore $A$ is included in the right hand side.

Suppose $n \ge 0$ and $(\alpha,\beta), (\alpha,\beta+n) \in S$.
We shall prove $v^\alpha w^\beta x^n \in A$ by induction on $n$.
It follows by definition when $n = 0, 1$.
Assume $n \ge 2$.
Then, by the induction hypothesis,
\[
v^{\alpha+1}w^\beta x^n - v^{\alpha}w^\beta x^n
= v^{\alpha}(v-1)w^\beta x^n = v^{\alpha}w^\beta(w-1) x^{n-1}
= v^{\alpha}w^{\beta+1} x^{n-1} - v^{\alpha}w^\beta x^{n-1} \in A
\]
since $(\alpha, \beta+1), (\alpha, \beta+n), (\alpha, \beta), (\alpha, \beta+n-1) \in S$. 
Therefore it is enough to show $v^{\alpha+1}w^\beta x^n\in A$.
Here remark $(\alpha+1, \beta), (\alpha+1, \beta+n) \in S$.
By this argument,  it is enough to show $v^{\alpha+q}w^\beta x^n\in A$
for $q \gg 0$.
Suppose $(q,1) \in S$.
Then $v^q x \in A$ since $(q,0), (q,1) \in S$.
We have
\[
v^{\alpha+q}w^\beta x^n = (v^q x)(v^{\alpha}w^\beta x^{n-1}) \in A .
\]

Next we shall prove 
\[
\sum_{\tiny \begin{array}{c} n \ge 0 \\ (\alpha,\beta), (\alpha,\beta+n) \in S \end{array}} Kv^\alpha w^\beta x^n = 
\bigoplus_{\alpha \ge 0} \bigoplus_{\tiny \begin{array}{c} n \ge 0 \\ (\alpha,\lceil \alpha \overline{u} \rceil+n) \in S \end{array}} Kx_{\alpha,n} .
\]
Here remark that $x_{\alpha,n}$'s are linearly independent by Proposition~\ref{strF}.
If $\alpha \ge 0$,  we have $(\alpha,\lceil \alpha \overline{u} \rceil) \in S$.
Therefore the right hand side is included in the left one.
Next we shall prove the opposite containment.
Suppose $n \ge 0$ and $(\alpha,\beta), (\alpha,\beta+n) \in S$.
We shall show that $v^\alpha w^\beta x^n$ is in the right hand side.
By Claim~\ref{claim}, we have
\[
v^\alpha w^\beta x^n \in \sum_{(\gamma, m) \in W_{\alpha,n,\beta-\alpha\overline{u}+n} \cap \bZ^2}Kx_{\gamma, m} .
\]
It is enough to show $\gamma \ge 0$, $m \ge 0$ and $(\gamma,\lceil \gamma \overline{u} \rceil+m) \in S$ for any $(\gamma, m) \in W_{\alpha,n,\beta-\alpha\overline{u}+n} \cap \bZ^2$.
Remember $\gamma \ge \alpha \ge 0$ and $m \ge n \ge 0$.
Put $q = \gamma - \alpha \ge 0$.
Since $(\gamma, m) \in W_{\alpha,n,\beta-\alpha\overline{u}+n}$, we have
\[
0 \le m \le \beta-\alpha\overline{u}+n -q \overline{u}
= \beta-\gamma\overline{u}+n .
\]
Adding $\lceil \gamma \overline{u} \rceil$, we have
\[
\lceil \gamma \overline{u} \rceil \le
\lceil \gamma \overline{u} \rceil + m \le
\beta+n + \left( \lceil \gamma \overline{u} \rceil -\gamma\overline{u} \right) .
\]
Since $0 \le \lceil \gamma \overline{u} \rceil -\gamma\overline{u} < 1$, we have
\[
\lceil \gamma \overline{u} \rceil \le
\lceil \gamma \overline{u} \rceil + m \le
\beta+n .
\]
Since $(\alpha,\beta+n) \in S$, we know $(\gamma,\beta+n) \in S$.
Therefore $(\gamma,\lceil \gamma \overline{u} \rceil+m) \in S$.
We have completed the proof of (\ref{strA1}).

(\ref{strA2}) follows from (\ref{strF1}), (\ref{strF2}) and (\ref{strA1}).
\qed

\vspace*{2mm}

For $\alpha \in \bZ$ and $n \ge 0$, we define
\[
y_{\alpha,n} := v^\alpha w^{\lceil \alpha \overline{u} \rceil} y^n \in G .
\]

\begin{Proposition}\label{strG}
We have 
\begin{align}
\label{strG1} G & = \bigoplus_{\alpha \in \bZ} \bigoplus_{n \ge 0} Ky_{\alpha,n} \\
\label{strG2}  y^\ell G & = \bigoplus_{\alpha \in \bZ} \bigoplus_{n \ge \ell } Ky_{\alpha,n} 
\end{align}
for any positive integer $\ell$.
\end{Proposition}

\proof
We have
\[
G = K[Z][y] = 
\sum_{(\alpha,\beta) \in Z\cap \bZ^2} \sum_{n\ge 0}K v^\alpha w^\beta y^n
\supset 
\sum_{\alpha \in \bZ} \sum_{n\ge 0}K y_{\alpha,n} .
\]

We shall prove the opposite containment.
For rational numbers $p$, $q$, $r$ such that $q<r$, let $W'_{p,q,r}$ be the triangle
with vertices ${\rm A}(p,q)$, ${\rm B}(p,r)$, ${\rm C}(\frac{r-q+p\overline{u}}{\overline{u}},\frac{q-r+q\overline{u}}{\overline{u}})$. 
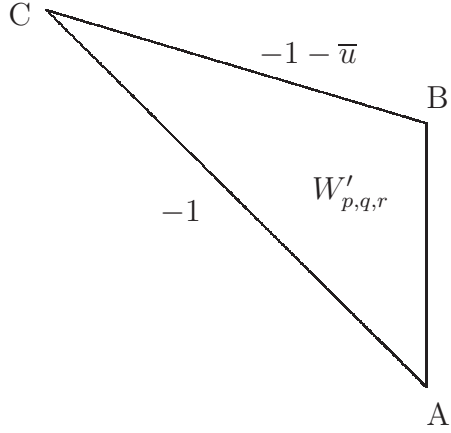
\[
{
\setlength\unitlength{1truecm}
  \begin{picture}(7,6)(-1,-2)
  %\put(-1,0){\vector(1,0){7}}
  %\put(0,-3){\vector(0,1){6}}
\qbezier (3,2) (0.5,2.75) (-2,3.5)
\qbezier (3,-1.5) (3,1) (3,2)
\qbezier (3,-1.5) (2,-0.5) (-2,3.5)
%\put(2.9,-1.6){$\bullet$}
\put(3,-2){A}
\put(3,2.2){B}
\put(-2.5,3.3){C}
  \put(-0.5,0.7){$-1$}
  \put(0.8,2.8){$-1-\overline{u}$}
  \put(1.5,1){$W'_{p,q,r}$}
      %\put(0.7,0.3){$T$}
  \end{picture}
}
\]
%\[
%\begin{tikzpicture}
% \coordinate[label=right:A] (A) at (0,0); 
% \coordinate[label=left:C] (C) at (-4,4); 
% \coordinate[label=right:B] (B) at (0,2.5); 
% \fill[lightgray] (A)--(B)--(C)--cycle; 
% \draw[very thick] (A)--(B)--(C)--cycle; 
%    \draw (-1.6,3.5)node{$-1-\overline{u}$};
%        \draw (-2.5,2)node{$-1$};
%  \draw (-1,1.8)node{$W'_{p,q,r}$};
%\end{tikzpicture}
%\]
The slope of the lower edge of $W'_{p,q,r}$ is $-1$, and that of the upper edge is $-1-\overline{u}$.

Let $\alpha$, $\beta$, $n$ be integers such that $\beta \ge \alpha \overline{u}$ and $n\ge 0$.
Then we can prove\footnote{
The coefficient of $y_{\alpha,n}$ is $1$.
If the characteristic of $K$ is $0$, we can prove
$v^\alpha w^\beta y^n \in \sum_{(\gamma, m) \in W'_{\alpha,n,\beta-\alpha\overline{u}+n} \cap \bZ^2}\bZ y_{\gamma, m}$.
} 
\[
v^\alpha w^\beta y^n \in \sum_{(\gamma, m) \in W'_{\alpha,n,\beta-\alpha\overline{u}+n} \cap \bZ^2}Ky_{\gamma, m} 
\]
in the same way as in Claim~\ref{claim}.

Thus (\ref{strG1}) and (\ref{strG2}) is proved as Proposition~\ref{strA}.
\qed

\begin{Remark}\label{algstrG}
\begin{rm}
For integers $\alpha$, $n$, let $V'_{\alpha,n}$ be the cone with vertex ${\rm A}(\alpha,n)$ as below:
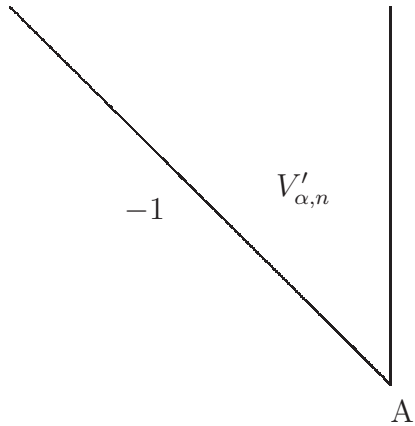
\[
{
\setlength\unitlength{1truecm}
  \begin{picture}(7,6)(-1,-2)
  %\put(-1,0){\vector(1,0){7}}
  %\put(0,-3){\vector(0,1){6}}
%\qbezier (3,2) (0.5,2.75) (-2,3.5)
\qbezier (3,-1.5) (3,1) (3,3.5)
\qbezier (3,-1.5) (2,-0.5) (-2,3.5)
%\put(2.9,-1.6){$\bullet$}
\put(3,-2){A}
%\put(3,2.2){B}
%\put(-2.5,3.3){C}
  \put(-0.5,0.7){$-1$}
  %\put(0.8,3){$-1-\overline{u}$}
  \put(1.5,1){$V'_{\alpha,n}$}
      %\put(0.7,0.3){$T$}
  \end{picture}
}
\]
%\[
%\begin{tikzpicture}
%\tikz \filldraw[fill=lightgray,very thick] (0,4)--(4,0)--(4,4);
% \coordinate[label=right:A] (A) at (4,0); 
%         \draw (3,2.2)node{$V'_{\alpha,n}$};
%  \draw (1.5,1.8)node{$-1$};
%\end{tikzpicture}
%\]

By the proof of Proposition~\ref{strG}, we have\footnote{
If the characteristic of $K$ is $0$, we can prove
$y_{\alpha,n}y_{\alpha',n'} = v^{\alpha+\alpha'}w^{\lceil \alpha \overline{u} \rceil + \lceil \alpha' \overline{u} \rceil}y^{n+n'}\in \sum_{(\gamma,m) \in V'_{\alpha+\alpha',n+n'} \cap \bZ^2}
\bZ y_{\gamma,m}$.
} 
\[
y_{\alpha,n}y_{\alpha',n'} = v^{\alpha+\alpha'}w^{\lceil \alpha \overline{u} \rceil + \lceil \alpha' \overline{u} \rceil}y^{n+n'}\in \sum_{(\gamma,m) \in V'_{\alpha+\alpha',n+n'} \cap \bZ^2}
Ky_{\gamma,m} ,
\]
where the coefficient of $y_{\alpha+\alpha',n+n'}$ is $1$.
Here remark $\lceil (\alpha  + \alpha' )\overline{u} \rceil \le \lceil \alpha \overline{u} \rceil + \lceil \alpha' \overline{u} \rceil$.
\end{rm}
\end{Remark}

\begin{Proposition}\label{strB}
We have the following equalities:
\begin{align}
\label{strB1}
B & = \sum_{\tiny \begin{array}{c} n \ge 0 \\ (\alpha,\beta), (\alpha,\beta+n) \in T \end{array}} Kv^\alpha w^\beta y^n = 
\bigoplus_{\alpha \le 0} \bigoplus_{\tiny \begin{array}{c} n \ge 0 \\ (\alpha,\lceil \alpha \overline{u} \rceil+n) \in T \end{array}} Ky_{\alpha,n} \\
\label{strB2}
x^\ell G \cap B & = \sum_{\tiny \begin{array}{c} n \ge \ell  \\ (\alpha,\beta), (\alpha,\beta+n) \in T \end{array}} Kv^\alpha w^\beta y^n= \bigoplus_{\alpha \le 0} \bigoplus_{\tiny \begin{array}{c} n \ge \ell  \\ (\alpha,\lceil \alpha \overline{u} \rceil+n) \in T \end{array}} Ky_{\alpha,n} 
\end{align}
for any positive integer $\ell$.
\end{Proposition}

We omit a proof since we can prove it in the same way as Proposition~\ref{strA}.

\vspace{2mm}

We put
\begin{align*}
F_\ell & = F/x^\ell F, \\
G_\ell & = G/y^\ell G, \\
A_\ell & = A/(x^\ell F \cap A), \\
B_\ell & = B/(y^\ell G \cap B).
\end{align*}
By (\ref{FtoG}), we have
\[%begin{equation}\label{BtoF}
B_\ell \subset G_\ell \simeq F_\ell \supset A_\ell .
\]%end{equation}
Let $$\psi:B_\ell \rightarrow F_\ell$$ be the above inclusion.
Then, since $-v^{-1}(1-x) = \frac{v^{-1}w-1}{v-1}$, we know
\[
\psi(y) = -\frac{vx}{1-x}
\]
and
\[
\psi(y_{\alpha,n}) = v^\alpha w^{\lceil \alpha \overline{u} \rceil} 
\left( -\frac{vx}{1-x} \right)^n
= (-1)^n v^{\alpha+n}w^{\lceil \alpha \overline{u} \rceil} (x+x^2+x^3+\cdots)^n .
\]
For $\alpha\in \bZ$ and $n \in \bNo$, we put
\[
z_{\alpha,n} := (-1)^n\psi(y_{\alpha-n,n}) = v^{\alpha}w^{\lceil (\alpha-n) \overline{u} \rceil} (x+x^2+x^3+\cdots)^n .
\]
Since $n \ge 0$, $\lceil (\alpha-n) \overline{u} \rceil \ge \lceil \alpha \overline{u} \rceil$.
By Claim~\ref{claim}, we can describe\footnote{
If the characteristic of $K$ is $0$, we can prove
$c_{\gamma,m} \in \bZ$.
} 
\begin{equation}\label{znosiki}
z_{\alpha,n} = \sum_{(\gamma,m) \in V_{\alpha,n}\cap \bZ^2}c_{\gamma,m}
x_{\gamma,m} \ \ (c_{\gamma,m} \in K)
\end{equation}
with $c_{\alpha,n} = 1$.
Then, we have
\begin{align}
F_\ell & = \bigoplus_{\alpha \in \bZ} \bigoplus_{\ell > n \ge 0} Kx_{\alpha,n} , \notag \\
\label{Anosiki} A_\ell & = \bigoplus_{\alpha \ge 0} \bigoplus_{\tiny
\begin{array}{c}
\ell > n \ge 0 \\ (\alpha, \lceil \alpha \overline{u} \rceil+n) \in S 
\end{array}
} Kx_{\alpha,n} , \\
\label{Bnosiki} \psi(B_\ell) & = 
\bigoplus_{\alpha \in \bZ} \bigoplus_{\tiny
\begin{array}{c}
\ell > n \ge 0 \\ (\alpha, \lceil \alpha \overline{u} \rceil+n) \in T
\end{array}
} Ky_{\alpha,n}
=\bigoplus_{\alpha \in \bZ} \bigoplus_{\tiny
\begin{array}{c}
\ell > n \ge 0 \\ (\alpha-n, \lceil (\alpha-n) \overline{u} \rceil+n) \in T
\end{array}
} Kz_{\alpha,n} .
\end{align}

We have
\[
(\spec A)|_{\ell C} = \spec A_\ell, \ \ 
(\spec B)|_{\ell C} = \spec B_\ell
\]
and
\[
\ell C = \spec A_\ell \cup \spec B_\ell, \ \ \spec F_\ell = \spec A_\ell \cap \spec B_\ell
\]
by (\ref{coverofellC}).

By Theorem~\ref{new}, it is very important whether $\oo_Y(m(abH-uE))|_{\ell C}$ is
a free sheaf or not.
We denote $\oo_Y(m(abH-uE))|_{\ell C}$ by $\lcal_{m, \ell}$.
As we have seen just before Theorem~\ref{new}, it is a locally free sheaf on $\ell C$.
Here put
\[
\xi := v^uw^{-u_2} \left( \frac{v^{-1}w-1}{1-v} \right)^u
= \left( \frac{w-v}{1-v} \right)^u w^{-u_2}
= (1-x)^u (1-x+vx)^{-u_2}
=(1-x_{0,1})^u(1-x_{0,1}+x_{1,1})^{-u_2}
.
\]
Remark that the constant term of $\xi$ is $1$.
We know that
\begin{equation}\label{batu}
\begin{array}{cl}
\bullet & \mbox{$\lcal_{m, \ell}|_{\spec A_\ell} \simeq \oo_{\spec A_\ell}$,} \\
\bullet & \mbox{$\lcal_{m, \ell}|_{\spec B_\ell} \simeq \oo_{\spec B_\ell}$,} \\
\bullet & \mbox{the transition function of $\lcal_{m, \ell}$ is
$\xi^m \in F_\ell^\times$.}
\end{array}
\end{equation}
Here $( \ )^\times$ is the set of unit elements in the given ring.
Thus we obtain the following:

\begin{Proposition}\label{criterion}
$\lcal_{m, \ell} \simeq \oo_{\ell C} \Longleftrightarrow
\mbox{$\exists \xi_A \in A_\ell^\times$, $\exists \xi_B \in B_\ell^\times$
such that $\xi^m = \xi_A\xi_B \in F_\ell^\times$}$
\end{Proposition}

\begin{Remark}\label{criteria}
\begin{rm}
Let $m$ be a positive integer.
Then the following conditions are equivalent:
\begin{enumerate}
\item
$\oo_Y(m(abH-uE))|_{mu C} \simeq \oo_{mu C}$.
\item
$\exists \xi_A \in A_{mu}^\times$, $\exists \xi_B \in B_{mu}^\times$
such that $\xi^m = \xi_A\xi_B \in F_{mu}^\times$.
\item
There exists $g \in [{\frak p}^{(mu)}]_{mab}$
such that $z^u-x^{s_3}y^{t_3}$, $g$ satisfying Huneke's criterion~\cite{Hu},
that is, $\ell_P(P/(z^u-x^{s_3}y^{t_3}, g, x)) = mua$ holds. 
\item
There exists an effective Weil divisor $D$ such that $D \sim m(abH-uE)$ and $C\cap D = \emptyset$.
\end{enumerate}
\end{rm}
\end{Remark}

\section{Strategy to prove finite/infinite generation}
\label{Strategy}

%Here we put
%\begin{align*}
%P_A & = \{ (\alpha,n) \in \bZ^2 \mid \alpha \ge 0, \ n \ge 0, \ (\alpha, \lceil \alpha \overline{u} \rceil+n) \in S \} , \\
%P_B & = \{ (\alpha,n) \in \bZ^2 \mid \alpha \in \bZ, \ n \ge 0, \ (\alpha-n, \lceil (\alpha-n) \overline{u} \rceil+n) \in T \} .
%\end{align*}

For $i = 1, 2, \ldots, mu$, we put
\[
q_i = ^\# \{ (\alpha, \beta) \in m\Delta_{\overline{t}, \overline{u}, \overline{s}} \cap {\Bbb Z}^2 \mid \alpha = i \} .
\]
Note that $q_{mu}=1$ and $q_i \ge 1$ for all $i = 1, 2, \ldots, mu$.
We sort the sequence $q_1$, $q_2$, \ldots, $q_{mu}$ into ascending order
\begin{equation}\label{q}
q'_1 \le q'_2 \le \cdots \le q'_{mu} .
\end{equation}

We say that the {\em condition EMU} is satisfied for $m\Delta_{\overline{t}, \overline{u}, \overline{s}}$ if
\[
q'_i \ge i
\]
for $i = 1, 2, \ldots, mu$.
If the condition EMU for $m\Delta_{\overline{t}, \overline{u}, \overline{s}}$ is satisfied
for some $m>0$, then $R_s({\frak p})$ is Noetherian.
One can prove it in the same way as the proof of Proposition~4.6 in \cite{KN}.

By definition, the condition EMU is satisfied for $(a,b,c)$ iff the condition EMU is satisfied for $\Delta_{\overline{t}, \overline{u}, \overline{s}}$.

Let $(a,b,c)$ be $(53,48,529)$.
It satisfies (\ref{abcjyouken}).
It is not difficult to prove that the condition EMU is satisfied for $(53,48,529)$.
Therefore $R_s({\frak p})$ is Noetherian.
However the condition EMU for $m\Delta_{\overline{t}, \overline{u}, \overline{s}}$ is not satisfied for any $m\ge 2$.

The following lemma is proved by Zhuang He.
We give a proof of it for a reader below.

\begin{Lemma}[He~\cite{He}]\label{He}
Assume that the condition EMU for $m\Delta_{\overline{t}, \overline{u}, \overline{s}}$ 
is not satisfied.
Then there exists an integer $d$ satisfying $1 \le d < mu$ such that the sequence
$q'_1$, $q'_2$, \ldots, $q'_{mu}$ defined in (\ref{q}) coincides with
\[
1, \ 2, \ 3, \ \ldots, \ d-1, \ d, \ d, \ldots .
\]
\end{Lemma}

\proof
For $i \in \bZ$, we put
\begin{equation}\label{defab}
\begin{split}
a_i & = ^\# \{ (\alpha, \beta) \in S \cap {\Bbb Z}^2 \mid \alpha = i \}  \\
b_i & = ^\# \{ (\alpha, \beta) \in T \cap {\Bbb Z}^2 \mid \alpha = i \} .
\end{split}
\end{equation}
By definition, we have
\[
\cdots \ge b_{-3} \ge b_{-2} \ge b_{-1} \ge 1 \le a_1\le a_2\le a_3 \le \cdots
\]
and
\[
\lim_{i \to \infty}a_i=\lim_{i \to \infty}b_{-i} = \infty .
\]
Here we remark $a_1>0$ and $b_{-1}>0$ since $\overline{t}<-1<\overline{u}<0<\overline{s}$.
We sort the sequence $1, a_1, b_{-1}, a_2, b_{-2}, a_3, b_{-3}, a_4, b_{-4}, \ldots$ into ascending order
\[
c_1 \le c_2 \le \cdots \le c_{i} \le \cdots .
\]
For $i>0$, let $e_i$ be the positive integer such that $c_{e_i} \le i < c_{e_i+1}$.
Then it is easy to see the following:
\begin{align}
\label{seq} \ \ &
\mbox{$c_i = q'_i$ for $i = 1, 2, \ldots, mu$, where $q'_i$ is defined in (\ref{q}).} \\
 \ \ &
\label{seqq}\mbox{The condition EMU for $m\Delta_{\overline{t}, \overline{u}, \overline{s}}$ is not satisfied if and only if there exists $i$} \\ & 
\mbox{satisfying $0<i<mu$ and $e_i > i$.} \notag
\end{align}

By definition, we have
\begin{align*}
a_i & = \lfloor i\overline{s} \rfloor - \lceil i\overline{u} \rceil + 1 \\
b_{-i} & = \lfloor -i\overline{t} \rfloor - \lceil -i\overline{u} \rceil + 1 
\end{align*}
for $i > 0$, where $\lfloor i\overline{s} \rfloor$ is the maximal integer such that $i\overline{s} \ge \lfloor i\overline{s} \rfloor$.
Since 
\begin{align*}
\lfloor (i+j)\overline{s} \rfloor & \ge \lfloor i\overline{s} \rfloor + \lfloor j\overline{s} \rfloor \\
\lceil (i+j)\overline{u} \rceil & \le \lceil i\overline{u} \rceil + \lceil j\overline{u} \rceil ,
\end{align*}
we obtain 
\begin{equation}\label{abplus}
\begin{array}{l}
a_{i+j} \ge a_i+a_j- 1 \\
b_{-(i+j)} \ge b_{-i}+b_{-j}-1 
\end{array}
\end{equation}
for positive integers $i$ and $j$.

If $a_1 = 1$ or $b_{-1} = 1$, then $q'_1 = q'_2 = 1$.
Here remark $u \ge 2$ and (\ref{seq}).

Next assume $a_1\ge 2$ and $b_{-1} \ge 2$.
By (\ref{abplus}), we have
\[
1 < a_1< a_2< a_3 < \cdots
\]
and
\[
1 < b_{-1}< b_{-2}< b_{-3} < \cdots .
\]
Since the condition EMU for $m\Delta_{\overline{t}, \overline{u}, \overline{s}}$ 
is not satisfied,
the set
\[
\{ a_1, a_2, a_3, \ldots \} \cap \{ b_{-1}, b_{-2}, b_{-3}, \ldots \}
\]
is not empty.
Let $d$ be the minimal number of the above set.
Since the condition EMU for $m\Delta_{\overline{t}, \overline{u}, \overline{s}}$ 
is not satisfied, we know $d<mu$.
Suppose $a_f = b_{-f'} = d$.
Then
\begin{equation}\label{inclusion}
\{ a_1, a_2, a_3, \ldots, a_{f-1} \} \coprod \{ b_{-1}, b_{-2}, b_{-3}, \ldots, b_{-(f'-1)} \}
\subset \{ 2, 3, \ldots, d-1 \} .
\end{equation}
By (\ref{seq}), it is enough to show that the above is the equality, that is, $f+f'-2 = d-2$.
 Assume the contrary, i.e. $d \ge f+f'+1$.

By (\ref{abplus}), we have $a_{nf} \ge nd - (n-1)$ and $b_{-nf'} \ge nd - (n-1)$.
Since $d \ge f+f'+1$, we have $nd-(n-1) \ge n(f+f')+1$.
Therefore we know 
\begin{equation}\label{abnokagen}
a_{nf} \ge n(f+f')+1,\ \ b_{-nf'} \ge n(f+f')+1 .
\end{equation}
Hence we have 
\[
e_{n(f+f')+1} \le n(f+f')+1, \ \ e_{n(f+f')} \le n(f+f')-1
\]
for any $n \ge 1$ by (\ref{seq}).
By (\ref{abplus}) and (\ref{abnokagen}), we have 
\[
a_{nf+1}\ge n(f+f')+a_1, \ a_{nf+2}\ge n(f+f')+a_2, \ \ldots, \ 
a_{nf+(f-1)}\ge n(f+f')+a_{f-1}
\]
and
\[
b_{-(nf'+1)}\ge n(f+f')+b_{-1}, \ b_{-(nf'+2)}\ge n(f+f')+b_{-2}, \ \ldots, \ 
b_{-(nf'+(f'-1))}\ge n(f+f')+b_{-(f'-1)} .
\]
Then, by (\ref{inclusion}), we know
\[
e_{n(f+f')+2} \le n(f+f')+2, \  e_{n(f+f')+3} \le n(f+f')+3, \ \ldots, \
e_{(n+1)(f+f')-1} \le (n+1)(f+f')-1 .
\]
It contradicts to (\ref{seqq}).
\qed

\begin{Definition}\label{mindeg}
\begin{rm}
Assume that the condition EMU for $m\Delta_{\overline{t}, \overline{u}, \overline{s}}$ 
is not satisfied.
The integer $d$ as in Lemma~\ref{He} is called the {\em minimal degree}
of $m\Delta_{\overline{t}, \overline{u}, \overline{s}}$.
\end{rm}
\end{Definition}

The above definition is due to He~\cite{He}.

For integers $p$, $q$, we put
\[
[p,q]_\bZ := \{ n \in \bZ \mid p \le n \le q \} , \ \
[p,q)_\bZ := \{ n \in \bZ \mid p \le n < q \} .
\]

We put
\begin{align*}
P_A & := \left\{ (\alpha,n) \in \bZ^2 \ \left| \  
\begin{array}{l}
\alpha \ge 0, \ n \ge 0, \\
(\alpha, \lceil \alpha \overline{u} \rceil+n) \in S 
\end{array}
 \right.  \right\} 
 =
 \left\{ (\alpha,n) \in \bZ^2 \ \left| \  
\begin{array}{l}
\alpha \ge 0, \ n \ge 0, \\
a_\alpha \ge n+1 
\end{array}
 \right.  \right\} 
, \\
P_B & := 
\left\{ (\alpha,n) \in \bZ^2 \ \left| \  
\begin{array}{l}
\alpha \in \bZ, \ n \ge 0, \\
(\alpha-n, \lceil (\alpha-n) \overline{u} \rceil+n) \in T 
\end{array}
 \right.  \right\} 
 =
 \left\{ (\alpha,n) \in \bZ^2 \ \left| \  
\begin{array}{l}
\alpha \in \bZ, \ n \ge 0, \\
b_{\alpha-n} \ge n+1 
\end{array}
 \right.  \right\} .
\end{align*}

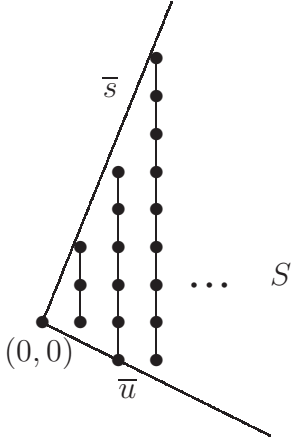
\[
{
\setlength\unitlength{1truecm}
  \begin{picture}(5.5,6)(2.5,-3)
  %\put(-1,0){\vector(1,0){7}}
  %\put(0,-3){\vector(0,1){6}}
\qbezier  (3,-1.5) (4,-2) (6,-3)
\qbezier (3,-1.5) (4,1) (4.7, 2.75)
\put(2.9,-1.6){$\bullet$}
\put(3.4,-1.6){$\bullet$}
\put(3.4,-1.1){$\bullet$}
\put(3.4,-0.6){$\bullet$}
\qbezier (3.5,-1.5) (3.5,-1) (3.5,-0.5)
\put(3.9,-2.1){$\bullet$}
\put(3.9,-1.6){$\bullet$}
\put(3.9,-1.1){$\bullet$}
\put(3.9,-0.6){$\bullet$}
\put(3.9,-0.1){$\bullet$}
\put(3.9,0.4){$\bullet$}
\qbezier  (4,-2) (4,-1.5) (4,0.5)
\put(4.4,-2.1){$\bullet$}
\put(4.4,-1.6){$\bullet$}
\put(4.4,-1.1){$\bullet$}
\put(4.4,-0.6){$\bullet$}
\put(4.4,-0.1){$\bullet$}
\put(4.4,0.4){$\bullet$}
\put(4.4,0.9){$\bullet$}
\put(4.4,1.4){$\bullet$}
\put(4.4,1.9){$\bullet$}
\qbezier (4.5,-2) (4.5,-1.5) (4.5,2)
\put(2.5,-2){$(0,0)$}
\put(6,-1){$S$}
  \put(4,-2.5){$\overline{u}$}
  \put(3.8,1.5){$\overline{s}$}
   \multiput(5,-1)(0.2,0){3}{\circle*{0.08}}
  \end{picture}
}
{
\setlength\unitlength{1truecm}
  
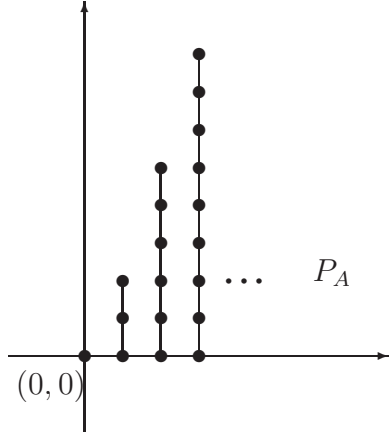
\begin{picture}(7,6)(-2,-3)
  \put(-1,-2){\vector(1,0){5}}
  \put(0,-3){\vector(0,1){5.7}}
\put(-0.1,-2.1){$\bullet$}
\put(0.4,-2.1){$\bullet$}
\put(0.4,-1.6){$\bullet$}
\put(0.4,-1.1){$\bullet$}
\qbezier (0.5,-2) (0.5,-1.5) (0.5,-1)
\put(0.9,-2.1){$\bullet$}
\put(0.9,-1.6){$\bullet$}
\put(0.9,-1.1){$\bullet$}
\put(0.9,-0.6){$\bullet$}
\put(0.9,-0.1){$\bullet$}
\put(0.9,0.4){$\bullet$}
\qbezier (1,-2) (1,-1.5) (1,0.5)
\put(1.4,-2.1){$\bullet$}
\put(1.4,-1.6){$\bullet$}
\put(1.4,-1.1){$\bullet$}
\put(1.4,-0.6){$\bullet$}
\put(1.4,-0.1){$\bullet$}
\put(1.4,0.4){$\bullet$}
\put(1.4,0.9){$\bullet$}
\put(1.4,1.4){$\bullet$}
\put(1.4,1.9){$\bullet$}
\put(-0.9,-2.5){$(0,0)$}
\qbezier (1.5,-2) (1.5,-1.5) (1.5,2)
   \multiput(1.9,-1)(0.2,0){3}{\circle*{0.08}}
   \put(3,-1){$P_A$}
  \end{picture}
}
\]

\[
{
\setlength\unitlength{1truecm}
  
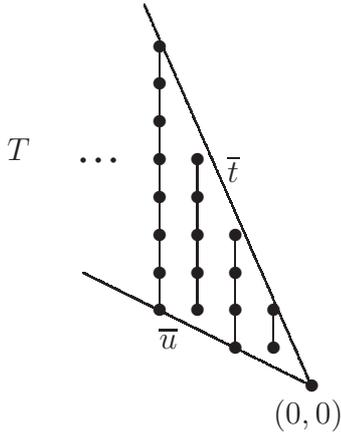
\begin{picture}(7,6)(-2,-2.5)
  %\put(-1,0){\vector(1,0){7}}
  %\put(0,-3){\vector(0,1){6}}
\qbezier (0,0) (1,-0.5) (3,-1.5) 
\qbezier (3,-1.5) (2,0.8) (0.8,3.56)
\put(2.9,-1.6){$\bullet$}
\put(2.4,-1.1){$\bullet$}
\put(2.4,-0.6){$\bullet$}
\qbezier (2.5,-1) (2.5,-0.7) (2.5,-0.5)
\put(1.9,-1.1){$\bullet$}
\put(1.9,-0.6){$\bullet$}
\put(1.9,-0.1){$\bullet$}
\put(1.9,0.4){$\bullet$}
\qbezier (2,-1) (2,0) (2,0.5)
\put(1.4,-0.6){$\bullet$}
\put(1.4,-0.1){$\bullet$}
\put(1.4,0.4){$\bullet$}
\put(1.4,0.9){$\bullet$}
\put(1.4,1.4){$\bullet$}
\qbezier (1.5,-0.5) (1.5,0) (1.5,1.5)
\put(0.9,-0.6){$\bullet$}
\put(0.9,-0.1){$\bullet$}
\put(0.9,0.4){$\bullet$}
\put(0.9,0.9){$\bullet$}
\put(0.9,1.4){$\bullet$}
\put(0.9,1.9){$\bullet$}
\put(0.9,2.4){$\bullet$}
\put(0.9,2.9){$\bullet$}
\qbezier (1,-0.5) (1,0) (1,3)
\put(2.5,-2){$(0,0)$}
  \put(1,-1){$\overline{u}$}
  \put(1.9,1.2){$\overline{t}$}
      \put(-1,1.5){$T$}
         \multiput(0,1.5)(0.2,0){3}{\circle*{0.08}}
  \end{picture}
}
{
\setlength\unitlength{1truecm}
  
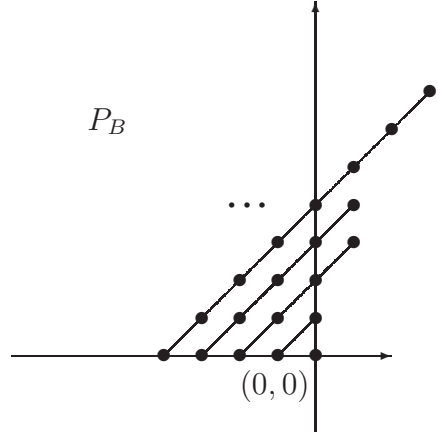
\begin{picture}(7,6)(-2,-3)
  \put(-1,-2){\vector(1,0){5}}
  \put(3,-3){\vector(0,1){5.7}}
\put(2.9,-2.1){$\bullet$}
\put(2.4,-2.1){$\bullet$}
\put(2.9,-1.6){$\bullet$}
\qbezier (2.5,-2) (2.75,-1.75) (3,-1.5)
\put(1.9,-2.1){$\bullet$}
\put(2.4,-1.6){$\bullet$}
\put(2.9,-1.1){$\bullet$}
\put(3.4,-0.6){$\bullet$}
\qbezier (2,-2) (2.5,-1.5) (3.5,-0.5)
\put(1.4,-2.1){$\bullet$}
\put(1.9,-1.6){$\bullet$}
\put(2.4,-1.1){$\bullet$}
\put(2.9,-0.6){$\bullet$}
\put(3.4,-0.1){$\bullet$}
\qbezier (1.5,-2) (2,-1.5) (3.5,0)
\put(0.9,-2.1){$\bullet$}
\put(1.4,-1.6){$\bullet$}
\put(1.9,-1.1){$\bullet$}
\put(2.4,-0.6){$\bullet$}
\put(2.9,-0.1){$\bullet$}
\put(3.4,0.4){$\bullet$}
\put(3.9,0.9){$\bullet$}
\put(4.4,1.4){$\bullet$}
\qbezier (1,-2) (1.5,-1.5) (4.5,1.5)
\put(2,-2.5){$(0,0)$}
 \multiput(1.9,0)(0.2,0){3}{\circle*{0.08}}
   \put(0,1){$P_B$}
  \end{picture}
}
\]

\begin{Lemma}\label{EMU}
\begin{enumerate}
\item
Let $m$ be a positive integer.
The following two conditions are equivalent:
\begin{enumerate}
\item
The condition EMU is satisfied for $m\Delta_{\overline{t}, \overline{u}, \overline{s}}$.
\item
$\bZ \times [0,mu)_\bZ \subset P_A\cup P_B$.
\end{enumerate}
\item
The set $\bNo^2 \setminus (P_A\cup P_B)$ is finite.
\item
Assume that the condition EMU for $m\Delta_{\overline{t}, \overline{u}, \overline{s}}$ 
is not satisfied.
Then there exist integers $d$ and $f$ such that
\begin{itemize}
\item
$mu > d \ge f \ge 0$,
\item
$(f,d) \not\in P_A \cup P_B$
\item
$P_A \cup P_B \cup \{ (f,d) \} \supset [0,f]_\bZ \times [0,d]_\bZ$,
\item
$P_A \cap P_B \cap ([0,f]_\bZ \times [0,d]_\bZ) = \{ (0,0) \}$.
\end{itemize}
\end{enumerate}
\end{Lemma}

\Proof
First we shall prove (1).
Remember $u \ge 2$ since ${\frak p}$ is not complete intersection.
Therefore the condition (a) implies $a_1 \ge 2$ and $b_{-1} \ge 2$.
On the other hand, assume the condition (b).
Since $(0,1) \not\in P_A$, we have $(0,1) \in P_B$ and $b_{-1}\ge 2$.
Since $(1,1) \not\in P_B$, we have $(1,1) \in P_A$ and $a_1 \ge 2$. 

From now on, we assume $a_1 \ge 2$ and $b_{-1} \ge 2$.
Then we know 
\begin{equation}\label{aandb}
\begin{array}{l}
2 \le a_1<a_2<a_3<\cdots , \\
2\le b_{-1}<b_{-2}<b_{-3}<\cdots , \\
\mbox{if $(\alpha,n) \in P_A$, then $(\alpha+1,n), (\alpha+1,n+1) \in P_A$,} \\
\mbox{if $(\alpha,n) \in P_B$, then $(\alpha-1,n), (\alpha,n+1) \in P_B$}
\end{array}
\end{equation}
by (\ref{abplus}).
We choose $q$ such that $a_q < mu\le a_{q+1}$.
Then, by (\ref{aandb}), the condition (b) is equivalent to
\[
(0,a_0), (1,a_1), (2,a_2), \ldots, (q,a_q) \in P_B ,
\]
where we put $a_0=b_0 = 1$.
It is also equivalent to 
\[
b_{i-a_i} \ge a_i + 1
\]
for $i = 0, 1, 2, \ldots, q$.
Then, the condition (b) is equivalent to
\begin{itemize}
\item
$b_{i-a_i} \ge a_i + 1$, $b_{i-a_i -1} \ge a_i + 2$, \ldots, 
$b_{(i+1)-a_{i+1}+1}\ge a_{i+1}-1$ for $i = 0, 1, 2, \ldots, q-1$, and
\item
$b_{q-a_q} \ge a_q + 1$, $b_{q-a_q -1} \ge a_q + 2$, \ldots, 
$b_{q-mu+2}\ge
mu-1$.
\end{itemize}
It is equivalent to the condition (a).

Next we shall prove (2).
For $\alpha\in \bZ$ and $n\in \bNo$, we know
\[
(\alpha,n)\in P_A \Longleftrightarrow (\alpha, \lceil \alpha \overline{u} \rceil+n) \in S
\Longleftrightarrow 
\left\{
\begin{array}{l}
\alpha \ge 0 \\
0 \le n \le \alpha \overline{s} - \lceil \alpha \overline{u} \rceil .
\end{array}
\right.
\]
Since $-\alpha \overline{u}\ge -\lceil \alpha \overline{u} \rceil > -\alpha \overline{u}-1$,
we have
\begin{equation}\label{PA}
\left\{ (\alpha,n) \ \left| 
\begin{array}{l}
\alpha \ge 0 \\
0 \le n \le \alpha (\overline{s} - \overline{u} ) -1
\end{array}
\right. \right\}
\subset P_A \subset
\left\{ (\alpha,n) \ \left| 
\begin{array}{l}
\alpha \ge 0 \\
0 \le n \le \alpha (\overline{s} - \overline{u} )
\end{array}
\right. \right\}
.
\end{equation}
Furthermore we have
\[
(\alpha,n)\in P_B \Longleftrightarrow 
\left\{
\begin{array}{l}
n \ge 0 \\
(\alpha-n, \lceil (\alpha-n) \overline{u} \rceil+n) \in T
\end{array}
\right.
\Longleftrightarrow 
\left\{
\begin{array}{l}
n \ge \alpha \\
0 \le n \le (\alpha-n) \overline{t} - \lceil (\alpha-n) \overline{u} \rceil .
\end{array}
\right.
\]
By (\ref{stu}),  we know $\overline{u}-\overline{t} > 1$.
By $-(\alpha-n) \overline{u}\ge -\lceil (\alpha-n) \overline{u} \rceil > -(\alpha-n) \overline{u}-1$, we have
\begin{equation}\label{PB}
\left\{ (\alpha,n) \ \left| 
\begin{array}{l}
n \ge \alpha, \ \ n \ge 0 \\
n \ge \frac{\overline{u}-\overline{t}}{\overline{u}-\overline{t}-1} \alpha + \frac{1}{\overline{u}-\overline{t}-1}
\end{array}
\right. \right\}
\subset P_B \subset
\left\{ (\alpha,n) \ \left| 
\begin{array}{l}
n \ge \alpha, \ \ n \ge 0 \\
n \ge \frac{\overline{u}-\overline{t}}{\overline{u}-\overline{t}-1} \alpha 
\end{array}
\right. \right\}
.
\end{equation}
By (\ref{stu}), we have
\[
1 < \frac{\overline{u}-\overline{t}}{\overline{u}-\overline{t}-1}
< \overline{s} - \overline{u} .
\]
Then (2) immediately follows from (\ref{PA}) and (\ref{PB}).

Next we shall prove (3).
If $b_{-1} = 1$, then $(f,d) = (0,1)$ satisfies our requirements.
If $b_{-1} \ge 2$ and $a_1 = 1$, then $(f,d) = (1,1)$ satisfies our requirements.
From now on, we assume $a_1 \ge 2$ and $b_{-1} \ge 2$.
Then (\ref{aandb}) is satisfied.
Let $d$ be the minimal degree defined in Definition~\ref{mindeg}.
Suppose $a_f = b_{-f'} = d$.
Then $ d = f+f'$ by Lemma~\ref{He}.
We have $(0,0) \in P_A \cap P_B$ and $(f,d) \not\in P_A \cup P_B$
since $a_f=d$ and $b_{f-d} = b_{-f'} = d$.

Since $(0,1)\in P_B \setminus P_A$, we have $(0,q)\in P_B \setminus P_A$
for any $q > 0$.
Since $(f,d)\not\in P_B$, we have $(f,q)\in P_A \setminus P_B$
for any $0\le q <d$.

Suppose $0<p<f$.
We have $(p,a_p-1)\in P_A$ and $(p,a_p)\not\in P_A$.
Since 
\[
\{ b_{-1}, b_{-2}, \ldots, b_{-(f'-1)} \} =
\{ 2, 3, \ldots, d-1 \} \setminus \{ a_1, a_2, \ldots, a_{f-1} \}  ,
\]
we have $a_p - 1 \ge b_{p-a_p+1}$ and $a_p + 1 \le b_{p-a_p}$ .
Therefore we have $(p,a_p-1)\not\in P_B$ and $(p,a_p)\in P_B$.
Ultimately, we obtain
$(p,q) \in P_B\setminus P_A$ if $q\ge a_p$, and
$(p,q) \in P_A\setminus P_B$ if $0 \le q< a_p$.
\qed

\vspace{2mm}

We define
\[
F'_\ell := \bigoplus_{\niretu{\ell > n \ge 0}{(\alpha,n) \in V_{0,0} \cap \bZ^2}}
Kx_{\alpha,n} \subset F_\ell = F/x^\ell F
=
\bigoplus_{\alpha \in \bZ} \bigoplus_{\ell > n \ge 0}
Kx_{\alpha,n} .
\]
By (\ref{algstr}), $F'_\ell$ is a subring of $F_\ell$.
Recall (\ref{Anosiki}).
We define
\[
A'_\ell := A_\ell \cap F'_\ell
= \bigoplus_{\sanretu{\ell > n \ge 0}{(\alpha,n) \in V_{0,0} \cap \bZ^2}{(\alpha, \lceil \alpha \overline{u} \rceil +n)\in S}}
Kx_{\alpha,n} .
\]
Recall (\ref{Bnosiki}).
We define
\[
B'_\ell := \psi(B_\ell) \cap F'_\ell
= \sum_{\sanretu{\ell > n \ge 0}{(\alpha,n) \in V_{0,0} \cap \bZ^2}{(\alpha-n, \lceil (\alpha-n) \overline{u} \rceil +n)\in T}}
Kz_{\alpha,n} .
\]
The above equality follows from (\ref{znosiki}) and (\ref{Bnosiki}).

\begin{Lemma}\label{subring}
Let $\ell$ be a positive integer.
Suppose $\eta = \eta_A \eta_B$ for $\eta\in {F'_\ell}^\times$, $\eta_A \in A_\ell^\times$ and $\eta_B \in B_\ell^\times$.

Then we have 
$\eta_A \in {A'_\ell}^\times$ and $\eta_B \in {B'_\ell}^\times$.
\end{Lemma}

\proof
It is easy to see that the constant terms (the coefficient of $x_{0,0}$) of
$\eta$, $\eta_A$ and $\eta_B$ are not $0$
by (\ref{algstr}), (\ref{PA}) and (\ref{PB}).

Since  $\eta_A \eta_B \in F'_\ell$, it follows that 
$\eta_A$ and $\eta_B$ are in $F'_\ell$ by (\ref{algstr}), (\ref{znosiki}), (\ref{PA}) and (\ref{PB}).
\qed

\begin{Definition}\label{Ualphan}
\begin{rm}
For $0 \le \alpha \le n$, we define
\[
U_{\alpha,n} = \{ (\gamma,m) \in V_{0,0} \cap \bZ^2 \mid 0 < m < n \}
\cup \{ (0,n), (1,n), \ldots, (\alpha-1,n) \} .
\]
We say that 
\begin{equation}\label{tejun}
\mbox{$\displaystyle 
	\zeta = \sum_{(\gamma,m)\in V_{0,0} \cap \bZ^2}h_{\gamma,m}x_{\gamma,m}$ 
	is $0$ in $U_{\alpha,n}$}
\end{equation}
 if $h_{\gamma,m} = 0$ for $(\gamma,m) \in U_{\alpha,n}$.
\end{rm}
\end{Definition}

\begin{Remark}\label{sennryaku}
\begin{rm}
Suppose that the constant term of $\zeta \in F'_\ell$ is $1$.
% for some
%$n \ge p^e$, $\eta_A \in {A'_{p^eu}}^\times$ and $\eta_B \in {B'_{p^eu}}^\times$.
Take $(\alpha,n)\in V_{0,0}\setminus \{ (0,0) \}$.
We assume that $(\alpha,n) \in P_A\cup P_B$.
We do the following procedure:
\begin{equation}\label{tejun2}
\begin{array}{cl}
%\bullet & 
%\mbox{Assume that $\zeta$ is zero in $U_{\alpha,n}$ .} \\
\bullet & 
\mbox{Let $c$ be the coefficient of $x_{\alpha,n}$ in $\zeta$. } \\
\bullet & 
\mbox{If $(\alpha,n) \in P_A$, then we multiply $1 - cx_{\alpha,n}$ to $\zeta$.} \\
\bullet & 
\mbox{If $(\alpha,n) \in P_B\setminus P_A$, then we multiply $1 - cz_{\alpha,n}$ to $\zeta$.}
\end{array}
\end{equation}
Then we know that the element $\zeta'$ obtained after the above procedure satisfies
\begin{itemize}
\item
$\zeta - \zeta'$ is $0$ in $U_{\alpha,n}$,
\item
the coefficient of $x_{\alpha,n}$ in $\zeta'$ is $0$.
\end{itemize}

Then we know the following:
\begin{enumerate}
\item
Suppose that $\eta \in {F'_\ell}^\times$ is $0$ in $U_{\alpha,n}$.
If $U_{\gamma,m}\setminus U_{\alpha,n} \subset P_A\cup P_B$, 
there exists $\eta_A \in {A'_{\ell}}^\times$ and $\eta_B \in {B'_{\ell}}^\times$
such that $\eta \eta_A \eta_B$ is $0$ in $U_{\gamma,m}$.
\item
Suppose that $\eta \in {F'_\ell}^\times$ in $0$ in $U_{\alpha,n}$.
Assume that $(\alpha',n)\in P_A$ for $\alpha < \alpha' \le n$.
We know that 
there exists $\eta_A \in {A'_{\ell}}^\times$
such that $\eta \eta_A = 1+ex_{\alpha,n}$ in ${F'_{n+1}}^\times$
for some $e \in K$.
\end{enumerate}
\end{rm}
\end{Remark}

Remember that $\xi$ is the transition function of $\lcal_{1, u}$ as in (\ref{batu}).
The following is a special case of Cutkosky's theorem~\cite{C}.

\begin{Proposition}\label{ch=p>0}
Assume ${\rm ch}(K) = p > 0$.

Then there exist $\xi_A \in {A'_{p^eu}}^\times$ and $\xi_B \in {B'_{p^eu}}^\times$
such that $\xi^{p^e} = \xi_A\xi_B$ in ${F'_{p^eu}}^\times$ for $e\gg 0$.

In particular $R_s({\frak p})$ is Noetherian.
\end{Proposition}

\proof
By Lemma~\ref{EMU} (2), we may assume that the second coordinate of each points
in $\bNo^2 \setminus (P_A\cup P_B)$ is less than $p^e$.
Then $\xi^{p^e}$ is $0$ in $U_{0,p^e}$.
By Remark~\ref{sennryaku} (1), we know that $\xi^{p^e} \eta_A \eta_B$ is $0$ in $U_{0,p^eu}$ for some $\eta_A \in {A'_{p^eu}}^\times$ and $\eta_B \in {B'_{p^eu}}^\times$.
Then we have $\xi^{p^e} = c \eta_A^{-1} \eta_B^{-1}$ in ${F'_{p^eu}}^\times$
for $c \in K^\times$, 
$\eta_A^{-1} \in {A'_{p^eu}}^\times$ and $\eta_B^{-1} \in {B'_{p^eu}}^\times$.
\qed

\begin{Lemma}\label{nthroot}
Let $K$ be a field of characteristic $0$.
Let $C$ be a local $K$-algebra with $\dim_KC < \infty$.
Let ${\frak m}$ be the maximal ideal of $C$.

Let $\varphi \in C$ is an element such that $\varphi \equiv 1 \mod {\frak m}$.
Let $n$ be a positive integer.

Then there exists the unique element $\varphi' \in C$ such that
$\varphi' \equiv 1 \mod {\frak m}$ and $\varphi'^n = \varphi$.
\end{Lemma}

We omit a proof.
Here we give another proof to Theorem~1.1 in \cite{KN}.

\begin{Proposition}\label{ch=0}
Assume ${\rm ch}(K) = 0$.
Let $m$ be a positive integer.
Assume that there exist $\xi_A \in {A'_{mu}}^\times$ and $\xi_B \in {B'_{mu}}^\times$
such that $\xi^m = \xi_A\xi_B$ in ${F'_{mu}}^\times$.

Then there exist $\xi_{A,1} \in {A'_{u}}^\times$ and $\xi_{B,1} \in {B'_{u}}^\times$
such that $\xi = \xi_{A,1}\xi_{B,1}$ in ${F'_u}^\times$.
\end{Proposition}

\proof
We may assume that the constant terms of $\xi_A$, $\xi_B$ are $1$.
By Lemma~\ref{nthroot},
there exist $\xi'_A \in {A'_{mu}}^\times$ and $\xi'_B \in {B'_{mu}}^\times$
 such that their constant terms are $1$ and
${\xi'_A}^m = \xi_A$, ${\xi'_B}^m = \xi_B$. 
Then we have $\xi = \xi'_A\xi'_B$ in ${F'_u}^\times$ by the uniqueness.
\qed

\begin{Lemma}\label{strategy}
Assume that the condition EMU for $m\Delta_{\overline{t}, \overline{u}, \overline{s}}$ 
is not satisfied.
Choose $d$ and $f$ satisfying Lemma~\ref{EMU} (3).
Suppose that there exist $\eta_A \in {A'_{d+1}}^\times$ and $\eta_B \in {B'_{d+1}}^\times$ such that $\xi^m \eta_A\eta_B=1+cx_{f,d}$ in ${F'_{d+1}}^\times$ with $c \in K^\times$.

Then there do not exist $\xi_A \in {A'_{mu}}^\times$ and $\xi_B \in {B'_{mu}}^\times$ 
satisfying $\xi^m = \xi_A \xi_B$ in ${F'_{mu}}^\times$.
\end{Lemma}

\proof
Assume the contrary.
Then we have
\[
\xi_A \xi_B\eta_A\eta_B = (\xi_A \eta_A)(\xi_B\eta_B) =1+cx_{f,d}
\]
in ${F'_{d+1}}^\times$.
Let $U_{\alpha,n}$ be the set defined in Definition~\ref{Ualphan}.
Assume that both $\xi_A \eta_A$ and $\xi_B\eta_B$ are $0$ in $U_{\alpha,n}$,
and the coefficient of $x_{\alpha,n}$ in either $\xi_A \eta_A$ or $\xi_B\eta_B$ is not zero
for some $\alpha$ and $n$.
Therefore we have $(\alpha,n) \in P_A \cup P_B$ and $(\alpha,n)\neq (f,d)$.
Since $1+cx_{f,d}$ is $0$ in $U_{\alpha,n}$, 
we know that one of the following is satisfied:
\begin{itemize}
\item
$n < d$
\item
$\alpha < f$ and $n = d$
\end{itemize}
Therefore we know that only one of $P_A$ and $P_B$ contains $(\alpha,n)$.
If the coefficient of $x_{\alpha,n}$ in $\xi_A \eta_A$ is not $0$, then $(\alpha,n) \in P_A$.
Then $(\alpha,n) \not\in P_B$ and the coefficient of $x_{\alpha,n}$ in $(\xi_A \eta_A)(\xi_B\eta_B)$ is not $0$.
It is a contradiction.
If the coefficient of $x_{\alpha,n}$ in $\xi_B \eta_B$ is not $0$, then $(\alpha,n) \in P_B$.
Then $(\alpha,n) \not\in P_A$ and the coefficient of $x_{\alpha,n}$ in $(\xi_A \eta_A)(\xi_B\eta_B)$ is not $0$.
It is a contradiction.
\qed

\begin{Lemma}\label{strategy2}
Assume that the condition EMU for $\Delta_{\overline{t}, \overline{u}, \overline{s}}$ 
is not satisfied.
Choose $d$ and $f$ satisfying Lemma~\ref{EMU} (3) with $m=1$.
Furthermore assume 
\begin{equation}\label{onlyone}
\left( \bZ \times [0,u)_\bZ \right) \setminus (P_A\cup P_B) = \{ (f,d) \} .
\end{equation}
Then the following conditions are equivalent:
\begin{enumerate}
\item
If $\xi \eta_A\eta_B=1+cx_{f,d}$ in ${F'_{d+1}}^\times$  is satisfied
for some $c \in K$, $\eta_A \in {A'_{d+1}}^\times$ and $\eta_B \in {B'_{d+1}}^\times$,
then $c \neq 0$.
\item
There do not exist $\xi_A \in {A'_{u}}^\times$ and $\xi_B \in {B'_{u}}^\times$ 
satisfying $\xi = \xi_A \xi_B$ in ${F'_{u}}^\times$.
\end{enumerate}
\end{Lemma}

\proof
By Remark~\ref{sennryaku}, we can prove that 
there exist $c \in K$, $\eta_A \in {A'_{d+1}}^\times$ and $\eta_B \in {B'_{d+1}}^\times$
satisfying $\xi \eta_A\eta_B=1+cx_{f,d}$ in ${F'_{d+1}}^\times$.
Then, by Lemma~\ref{strategy}, we know $(1)\Rightarrow (2)$.

Next we shall prove the converse.
Assume that (1) is not satisfied.
Then
there exist $\eta_A \in {A'_{u}}^\times$ and $\eta_B \in {B'_{u}}^\times$ such that $\xi \eta_A\eta_B = 1$ in ${F'_{d+1}}^\times$.
By Remark~\ref{sennryaku} (1), we know that
there exist $\eta'_A \in {A'_{u}}^\times$ and $\eta'_B \in {B'_{u}}^\times$ such that $\xi \eta_A\eta_B\eta'_A\eta'_B=1$ in ${F'_{u}}^\times$ by the assumption (\ref{onlyone}).
Therefore we have $\xi =( \eta_A\eta'_A)^{-1}(\eta_B\eta'_B)^{-1}$ in ${F'_{u}}^\times$.
\qed

\section{Classification of $(a,b,c)'s$ which do not satisfy the condition EMU}

We want to classify $(a,b,c)$'s which do not satisfy the condition EMU.

\begin{Remark}\label{rem6.1}
\begin{rm}
If the condition EMU for $m\Delta_{\overline{t}, \overline{u}, \overline{s}}$ is satisfied
for some $m>0$, then $R_s({\frak p})$ is Noetherian.
One can prove it in the same way as the proof of Proposition~4.6 in \cite{KN}.

Remark that we may assume $\ell_1 \ge \ell_{u-1}$ by exchanging $x$ for $y$
if necessary. 

If $\ell_1=1$ or $\ell_{u-1} = 1$, the condition EMU for $\Delta_{\overline{t}, \overline{u}, \overline{s}}$ is not satisfied.
In this case, $R_s({\frak p})$ is not Noetherian by Theorem~1.2 in Gonz\'alez-Karu~\cite{GK}.

If $\ell_1=\ell_{u-1}=2$ and $u\ge 3$, then the condition EMU for $\Delta_{\overline{t}, \overline{u}, \overline{s}}$ is not satisfied.
In this case, $R_s({\frak p})$ is not Noetherian by Theorem~1.2 in Gonz\'alez-Karu~\cite{GK}.

If $\ell_1\ge 3$ and $\ell_{u-1} \ge 3$, then the condition EMU for $\Delta_{\overline{t}, \overline{u}, \overline{s}}$ is satisfied by Lemma~\ref{He}.

The remaining case is $\ell_1 = n \ge 3$ and $\ell_{u-1} = 2$.
\end{rm}
\end{Remark}

In the rest of this section, we suppose that $n$ is a positive integer such that $$n\ge 3 .$$
We want to classify $(a,b,c)$'s which do not satisfy the condition EMU for $(a,b,c)$ such that $a_1 = n$ and $b_{-1} = 2$.

First we define the sequence of integers
\[
\{ f_i, g_i \mid i = 0, 1, \ldots \}
\]
as
\[
\left(
\begin{array}{c}
f_0 \\ g_0
\end{array}
\right) 
=
\left(
\begin{array}{c}
1 \\ 1
\end{array}
\right) , \ \ 
\left(
\begin{array}{c}
f_i \\ g_i
\end{array}
\right) 
=
\left(
\begin{array}{cc}
n-2 & 1 \\ n-3 & 1
\end{array}
\right)^i
\left(
\begin{array}{c}
f_0 \\ g_0
\end{array}
\right) .
\]
We know 
\[
(f_0,g_0)=(1,1), \ (f_1,g_1)=(n-1,n-2), \  (f_2,g_2)=(n^2-2n,n^2-3n+1), \ \ldots .
\]
Here remark that 
\[
1 = f_0 < f_1 < f_2 < \cdots
\]
and
\[
1 = g_0 \le g_1 \le g_2 \le \cdots .
\]
For $i > 0$, we have
\[
\left(
\begin{array}{cc}
f_i & f_{i-1} \\ g_i & g_{i-1}
\end{array}
\right)
= 
\left(
\begin{array}{cc}
n-2 & 1 \\ n-3 & 1
\end{array}
\right)^{i-1}
\left(
\begin{array}{cc}
f_1 & f_0 \\ g_1& g_0
\end{array}
\right)
=
\left(
\begin{array}{cc}
n-2 & 1 \\ n-3 & 1
\end{array}
\right)^{i-1}
\left(
\begin{array}{cc}
n-1 & 1 \\ n-2& 1
\end{array}
\right) 
\]
and 
\begin{equation}\label{det}
{\rm det}\left(
\begin{array}{cc}
f_i & f_{i-1} \\ g_i & g_{i-1}
\end{array}
\right) = 1  .
\end{equation}
In particular, we have ${\rm GCD}(f_i,f_{i-1}) = {\rm GCD}(f_i,g_i) = 1$.

Since $f_i = (n-2)f_{i-1}+g_{i-1}$ and $g_i = (n-3)f_{i-1}+g_{i-1}$,
we know 
\begin{equation}\label{f=f-g}
f_i = g_i + f_{i-1}
\end{equation}
for $i > 0$.

For $i \ge 2$. we know
\begin{equation}\label{f=f-f}
f_i =  (n-2)f_{i-1}+g_{i-1} =  (n-2)f_{i-1}+f_{i-1}-f_{i-2} = (n-1)f_{i-1}-f_{i-2} .
\end{equation}

We define $a_i$ and $b_i$ as in (\ref{defab}).
By Lemma~\ref{He} and (\ref{seq}), the condition EMU is not satisfied for $\Delta_{\overline{t}, \overline{u}, \overline{s}}$ if and only if there exist positive integers $f$, $f'$, $d$ such that
\begin{itemize}
\item
$a_f = b_{-f'} = d$,
\item
$\{ a_1, a_2, \ldots, a_{f-1} \} \coprod \{ b_{-1}, b_{-2}, \ldots, b_{-(f'-1)} \}
= \{ 2, 3, \ldots, d-1 \}$,
\item
$f+f' =d<u$.
\end{itemize}

\begin{Theorem}\label{classification}
Let $a$, $b$, $c$ be pairwise coprime positive integers.
Assume that $\frak p$ is not complete intersection.
Let $n$ be an integer such that $n \ge 3$.
Then the following two conditions are equivalent:  

\begin{enumerate}
\item
The following conditions are satisfied:
\begin{enumerate}
\item
$z^u-x^{s_3}y^{t_3}$ is the negative curve, i.e., $uc < \sqrt{abc}$,
\item
the condition EMU is not satisfied for $\Delta_{\overline{t}, \overline{u}, \overline{s}}$,
\item
$\ell_1 = n$ and $\ell_{u-1} = 2$.
\end{enumerate}
\item
There exist $\lambda \in \bNo$ and $\gamma, \delta \in \bN$ with $(\gamma, \delta) \neq (1,1)$ and ${\rm GCD}(\gamma,\delta) = 1$ satisfying the following conditions\footnote{The minimal degree $d$ defined in Definition~\ref{mindeg} is $f_\lambda+f_{\lambda+1}$ in this case.}:
\begin{enumerate}
\item
$u = \gamma f_\lambda + \delta f_{\lambda+1}$, $u_2 = \gamma g_\lambda + \delta g_{\lambda+1}$ ,
\item
$n-1 \le \overline{s} < \frac{(n-1)f_\lambda+1}{f_\lambda}$,
\item
$2 \le -\overline{t} < \frac{2f_{\lambda+1}+1}{f_{\lambda+1}}$,
\item
$\frac{1}{\overline{s}-\overline{u}} + \frac{1}{\overline{u}-\overline{t}}<1$ .
\end{enumerate}
\end{enumerate}
\end{Theorem}

We shall prove this theorem in this section.

\begin{Definition}\label{gamma}
\begin{rm}
Let $f$, $p_1$, $p_2$, $f'$, $q_1$, $q_2$ be rational numbers such that $f>0$ and $f'>0$.
Let $\Gamma(-f',p_1,p_2;f,q_1,q_2)$ the set of lattice points consisting of
\[
\{ (\alpha,\beta) \in \bZ^2 \mid \alpha \le 0, \ \frac{-p_2\alpha}{f'} \le \beta \le \frac{-p_1\alpha}{f'} \} 
\]
and
\[
\{ (\alpha,\beta) \in \bZ^2 \mid 0<\alpha, \ \frac{q_2\alpha}{f} < \beta \le \frac{q_1\alpha}{f} \} .
\]
%Let $r_1$ and $r_2$ be positive integers.
%We say that $\Delta_{\overline{t}, \overline{u}, \overline{s}}$ coincides with 
%$\Gamma(f,p_1,p_2;f',q_1,q_2)$ in the range of $-r_2$ to $r_1$
%if the following two conditions are satisfied:
%\begin{itemize}
%\item
%$S \cap [0,r_1]_\bZ = \Gamma(f,p_1,p_2;f',q_1,q_2) \cap [0,r_1]_\bZ$.
%\item
%$T \cap [-r_2,0]_\bZ = \Gamma(f,p_1,p_2;f',q_1,q_2) \cap [-r_2,0]_\bZ$.
%\end{itemize}
\end{rm}
\end{Definition}

\begin{Lemma}\label{tuduku}
Let $f$ and $f'$ be positive integers.
Let $p_1$, $p_2$, $q_1$, $q_2$ be rational numbers.
For $i \in \bZ$, we define
\[
e_i = {}^\# \{ (i,\beta)\in \bZ^2 \mid (i,\beta) \in \Gamma(-f',p_1,p_2;f,q_1,q_2) \} .
\]
Assume that if we sort the sequence $e_{-f'}, e_{-f'+1}, \ldots, e_f$ into ascending order, we obtain 
\[
1, 2, \ldots, f+f', f+f'+1 .
\]
with $e_f = f + f'$ and $e_{-f'} = f+f'+1$.
Let $c_1, c_2, c_3, \ldots$ be the sequence given by sorting the sequence $\{ e_i \mid i \in \bZ \}$ into ascending order.
\begin{enumerate}
\item
Assume that $p_1$, $p_2$, $q_1$, $q_2$ are integers.
Then we have $c_i = i$ for any $i>0$.
\item
Assume that $p_1$, $q_1$, $q_2$ are integers and $p_2 \not\in \bZ$.
Then there exists $i > 0$ such that
\[
c_1=1, \ c_2 = 2, \ \ldots, \ c_{i-1} = i-1, \ c_i > i .
\]
\end{enumerate}
\end{Lemma}

\proof
First we prove (1).
We have $e_{f+i} = e_f + e_i = f+f'+e_i$ for $i > 0$ and
$e_{-(f'+i)} = e_{-f'}-1 + e_{-i} = f+f'+e_{-i}$ for $i \ge 0$.
Then we have $e_{2f} = 2(f + f')$ and $e_{-2f'} = 2(f+f')+1$.
Recall 
\[
\{ e_{-f'}, e_{-f'+1}, \ldots, e_f \} = \{ 1, 2, \ldots, f+f'+1 \} .
\]
Therefore, if we sort the sequence $e_{-2f'}, e_{-2f'+1}, \ldots, e_{2f}$ into ascending order, we obtain 
\[
1, 2, \ldots, 2(f+f'), 2(f+f')+1 .
\]
Repeating this process, we shall obtain the assertion.

Next we shall prove (2).
We have $e_{f+i} = e_f + e_i = f+f'+e_i$ for $i = 1, 2, \ldots, f$ and
\begin{equation}\label{ineq}
e_{-(f'+i)} \ge e_{-f'}-1 + e_{-i} = f+f'+e_{-i}
\end{equation}
for $i = 1, 2, \ldots, f'$.
If there exists $i$ such that $e_{-(f'+i)} > f+f'+e_{-i}$,
the assertion follows immediately.
If (\ref{ineq}) are the equalities for $i = 1, 2, \ldots, f'$, 
we have $e_{2f} = 2(f + f')$ and $e_{-2f'} = 2(f+f')+1$.
If we sort the sequence $e_{-2f'}, e_{-2f'+1}, \ldots, e_{2f}$ into ascending order, we obtain 
\[
1, 2, \ldots, 2(f+f'), 2(f+f')+1 .
\]
We repeat this process.
The equalities will not hold in the future.
\qed

\begin{Lemma}\label{tuduki}
Let $\lambda$ be a non-negative integer.
We put $f=f_\lambda$, $f' = f_{\lambda+1}$, $p_1 = 2f_{\lambda+1}$, $p_2 = g_{\lambda+1}$,
$q_1 = (n-1)f_\lambda$, $q_2 = -g_\lambda$ and
$\Gamma = \Gamma(-f_{\lambda+1}, 2f_{\lambda+1}, g_{\lambda+1}; f_{\lambda}, (n-1)f_{\lambda}, -g_{\lambda}) \}$.

\begin{enumerate}
\item
The integers $f$, $f'$, $p_1$, $p_2$,
$q_1$, $q_2$ satisfy the assumption of Lemma~\ref{tuduku}.
\item
The set of lattice points in $\Gamma$
with the first component $f_\lambda$ is
$$\{ (f_\lambda, -g_\lambda+1), (f_\lambda, -g_\lambda+2), \ldots, (f_\lambda, (n-1)f_\lambda) \} .$$
The number of this set is $f_\lambda + f_{\lambda+1}$.
\item
The set of lattice points in $\Gamma$
with the first component $-f_{\lambda+1}$ is
$$\{ (-f_{\lambda+1}, g_{\lambda+1}), (-f_{\lambda+1}, g_{\lambda+1}+1), \ldots, (-f_{\lambda+1}, 2f_{\lambda+1}) \} .$$
The number of this set is $f_\lambda + f_{\lambda+1}+1$.
\item
The set of lattice points in $\Gamma$
with the first component $f_{\lambda+1}$ is
$$\{ (f_{\lambda+1}, -g_{\lambda+1}), (f_{\lambda+1}, -g_{\lambda+1}+1), \ldots, (f_{\lambda+1}, (n-1)f_{\lambda+1}) \} .$$
The number of this set is $f_{\lambda+1} + f_{\lambda+2}+1$.
\item
The set of lattice points in $\Gamma$
with the first component $-f_{\lambda+2}$ is
$$\{ (-f_{\lambda+2}, g_{\lambda+2}+1), (-f_{\lambda+2}, g_{\lambda+2}+2), \ldots, (-f_{\lambda+2}, 2f_{\lambda+2}) \} .$$
The number of this set is $f_{\lambda+1} + f_{\lambda+2}$.
\end{enumerate}
\end{Lemma}

\proof
(2) and (3) follow from definition immediately.

We shall prove (1) by induction on $\lambda \ge 0$.

Assume $\lambda = 0$.
We have $f'=n-1$, $p_1=2(n-1)$, $p_2 = n-2$, $f=1$, $q_1=n-1$ and $q_2 = -1$.
Let $e_i$ be the number of lattice points in $\Gamma(-f',p_1,p_2;f,q_1,q_2)$
such that the first component is $i$.
Then we obtain 
$e_{-n+1}=n+1$, $e_{-n+2}=n-1$, $e_{-n+3}=n-2$, \ldots, $e_{-1}=2$, $e_0 = 1$, $e_1=n$.
Thus the assumption in Lemma~~\ref{tuduku} is satisfied.

Assume that the assumption in Lemma~~\ref{tuduku} is satisfied for some $\lambda\ge 0$.
We define
\begin{align*}
e'_i & = {}^\# \{ (i,\beta)\in \bZ^2 \mid 
(i,\beta)\in \Gamma(-f_{\lambda+1}, 2f_{\lambda+1}, g_{\lambda+1}; f_{\lambda}, (n-1)f_{\lambda}, -g_{\lambda}) \} , \\
e''_i & = {}^\# \{ (i,\beta)\in \bZ^2 \mid 
(i,\beta)\in \Gamma(-f_{\lambda+2}, 2f_{\lambda+2}, g_{\lambda+2}; f_{\lambda+1}, (n-1)f_{\lambda+1}, -g_{\lambda+1}) \} .
\end{align*}
Remark
\[
1 \le e'_1 \le e'_2 \le \cdots 
\]
and
\[
1 \le e'_{-1} \le e'_{-2} \le \cdots .
\]
By Lemma~\ref{tuduku} (1), 
if we sort the sequence $\{ e'_i \mid i \in \bZ\}$ into ascending order,
we obtain 
$$
1, \ 2, \ 3, \ \ldots .$$
Here we have
\begin{equation}\label{a'nosiki}
e''_{f_{\lambda+1}} = (n-1)f_{\lambda+1} + g_{\lambda+1} 
= (n-1)f_{\lambda+1} + (f_{\lambda+1}-f_\lambda)
= f_{\lambda+1}+ f_{\lambda+2}
\end{equation}
by (\ref{f=f-g}) and (\ref{f=f-f}).
Furthermore, we have
\begin{equation}\label{b'nosiki}
e''_{-f_{\lambda+2}} = 2f_{\lambda+2} - g_{\lambda+2}+1
= 2f_{\lambda+2} - (f_{\lambda+2}-f_{\lambda+1})+1
= f_{\lambda+1}+ f_{\lambda+2}+1 .
\end{equation}
By (\ref{det}), we have
\[
1= {\rm det}\left(
\begin{array}{cc}
f_{\lambda+1} & f_\lambda \\ g_{\lambda+1} & g_\lambda
\end{array}
\right) = 
{\rm det}\left(
\begin{array}{cc}
f_\lambda & f_{\lambda+1} \\ -g_\lambda & -g_{\lambda+1}
\end{array}
\right) .
\]
Therefore we know 
\[
e'_{f_{\lambda+1}} =e''_{f_{\lambda+1}} +1 = f_{\lambda+1}+ f_{\lambda+2}+1 
\]
and $e''_i = e'_i$ for $0<i < f_{\lambda+1}$.
Thus (4) follows from the above equality.
By (\ref{det}), we have
\[
1= {\rm det}\left(
\begin{array}{cc}
f_{\lambda+2} & f_{\lambda+1} \\ g_{\lambda+2} & g_{\lambda+1}
\end{array}
\right) = 
{\rm det}\left(
\begin{array}{cc}
-f_{\lambda+1} & -f_{\lambda+2} \\ g_{\lambda+1} & g_{\lambda+2}
\end{array}
\right) .
\]
Therefore we know
\[
e'_{-f_{\lambda+2}} =e''_{-f_{\lambda+2}} -1 = f_{\lambda+1}+ f_{\lambda+2} 
\]
and
 $e''_{i} = e'_i$ for $-f_{\lambda+2}<i<0$.
Thus (5) follows from the above equality. 
If we sort the sequence $e''_{-f_{\lambda+2}+1}, e''_{-f_{\lambda+2}+2}, \ldots, e''_{f_{\lambda+1}-1}$ into ascending order, we have
\[
1,2, \ldots, f_{\lambda+1}+ f_{\lambda+2}-1.
\]
By (\ref{a'nosiki}) and (\ref{b'nosiki}), we know that (1) holds for $\lambda + 1$.
\qed

\vspace{2mm}

Now we start to prove $(1) \Rightarrow (2)$ in Theorem~\ref{classification}.
We put $$f_{-1} = 0.$$
It is enough to prove the following claim:

\begin{Claim}\label{classify}
Assume that (1) in Theorem~\ref{classification} is satisfied.
Let $\lambda$ be non-negative integer.
Then we have the following:
\begin{itemize}
\item[$(1)_\lambda$]
Letting $d$ be the minimal degree of $\Delta_{\overline{t}, \overline{u}, \overline{s}}$,
$d$ does not satisfy $f_{\lambda-1}+f_{\lambda} < d < f_{\lambda}+f_{\lambda+1}$.
\item[$(2)_\lambda$]
If the minimal degree of $\Delta_{\overline{t}, \overline{u}, \overline{s}}$ is $f_{\lambda}+f_{\lambda+1}$, then 
there exist $\gamma, \delta \in \bN$ with $(\gamma, \delta) \neq (1,1)$ and ${\rm GCD}(\gamma,\delta) = 1$ such that 
$u$, $u_2$, $\overline{t}$,
$\overline{s}$ satisfy (a), (b), (c), (d) in (2) in Theorem~\ref{classification}.
\item[$(3)_\lambda$]
If the minimal degree of $\Delta_{\overline{t}, \overline{u}, \overline{s}}$ is bigger than
$f_{\lambda}+f_{\lambda+1}$, 
then  
%$\Delta_{\overline{t}, \overline{u}, \overline{s}}$
%coincides with $\Gamma(f_{\lambda},(n-1)f_{\lambda},-g_{\lambda};
%f_{\lambda+1}, 2f_{\lambda+1}, g_{\lambda+1} )$ in the range of $-f_{\lambda+1}$ to $f_{\lambda}$.
we have
\[
(S\cup T) \cap [-f_{\lambda+1},f_{\lambda}]_\bZ^1
= \Gamma(-f_{\lambda+1}, 2f_{\lambda+1}, g_{\lambda+1};f_{\lambda},(n-1)f_{\lambda},-g_{\lambda})
\cap [-f_{\lambda+1},f_{\lambda}]_\bZ^1 ,
\]
where $[-f_{\lambda+1},f_{\lambda}]_\bZ^1:=[-f_{\lambda+1},f_{\lambda}]_\bZ \times \bZ = \{ (\alpha,\beta) \in \bZ^2 \mid -f_{\lambda+1}\le \alpha \le f_{\lambda} \}$.
\end{itemize}
\end{Claim}

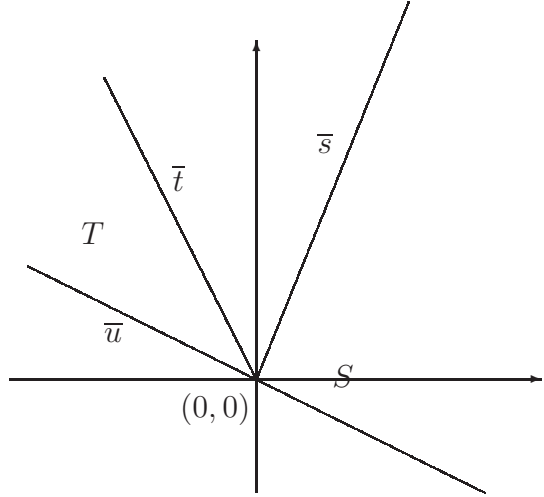
\[
{
\setlength\unitlength{1truecm}
  \begin{picture}(7,6)(0,-3)
  \put(-0.25,-1.5){\vector(1,0){7}}
  \put(3,-3){\vector(0,1){6}}
\qbezier (0,0) (3,-1.5) (6,-3)
\qbezier (3,-1.5) (2,0.5) (1,2.5)
\qbezier (3,-1.5) (4,1) (5, 3.5)
%\put(2.9,-1.6){$\bullet$}
\put(2,-2){$(0,0)$}
\put(4,-1.6){$S$}
  \put(1,-1){$\overline{u}$}
  \put(1.9,1){$\overline{t}$}
  \put(3.8,1.5){$\overline{s}$}
      \put(0.7,0.3){$T$}
  \end{picture}
}
\]
%\[
%\begin{tikzpicture}[x=2cm,y=2cm]
% \filldraw[fill=lightgray,very thick] (-1.5,0.5)--(0,0)--(-1.5,2.5);
% \filldraw[fill=lightgray,very thick] (1.5,-0.5)--(0,0)--(1.5,3);   
%\draw[->,>=stealth,semithick] (-1.7,0)--(1.7,0)node[below]{}; %x軸
% \draw[->,>=stealth,semithick] (0,-0.5)--(0,3)node[right]{}; %y軸
% \draw (0,0)node[below left]{$(0,0)$}; %原点 
%    \coordinate (A) at (0.3,1) node at (A) [above right] {$\overline{s}$};
%   \coordinate (P) at (0.7,0.5) node at (P) [below=0] {$S$};
%   \coordinate (B) at (-0.7,0.7) node at (B) [below=0] {$T$};
%   \coordinate (C) at (-0.3,1) node at (C) [below=0] {$\overline{t}$};
%   \coordinate (D) at (0.7,-0.3) node at (D) [below=0] {$\overline{u}$};
%%   \fill (C) circle [radius=2pt];
%% \draw[black,very thick,domain=-1.3:1.3] plot(\x,{(\x)^2 + (1/2)})node[left]{};
%% \draw[black,very thick,domain=-1.5:0.88] plot(\x,{1/(1-\x)})node[right]{$y=1/(1-x)$};
%%  \draw[red,very thick,domain=-1.1:1.1] plot(\x,{1/2})node[right]{};
%%  \draw[red,very thick,domain=-0.2:1.1] plot(\x,{(\x) + (1/4)})node[right]{};
%%    \draw[red,very thick,domain=-1.1:0.2] plot(\x,{(-1)*(\x) + (1/4)})node[right]{};
%%  \draw[red,very thick,domain=-3:3] plot(\x,{e^((-1)*(\x)^2)+0.05})node[above]{$y=e^{-x^2}$};
%\end{tikzpicture}
%\]

We shall prove this claim by induction on $\lambda$.

First, assume $\lambda = 0$.
Then we know $f_{-1}+f_0 = 1$ and $f_0+f_1 = n$.

Since $a_1 = n$, the minimal degree of $\Delta_{\overline{t}, \overline{u}, \overline{s}}$ 
is not less than $n$ by (\ref{aandb}) and Lemma~\ref{He}.
Thus $(1)_0$ follows.

Since $(-1,1), (-1,2) \in T$, we know 
%$\Delta_{\overline{t}, \overline{u}, \overline{s}}$
%coincides with $\Gamma(1,n-1,-1;n-1, 2(n-1), n-1)$ in the range of $-(n-1)$ to $1$.
\begin{equation}\label{e=0}
(S \cup T)\cap [-(n-1),1]_\bZ^1 \supset \Gamma(-1, 2, 1;1,n-1,-1)\cap  [-(n-1),1]_\bZ^1 .
\end{equation}
Assume that the minimal degree of $\Delta_{\overline{t}, \overline{u}, \overline{s}}$ 
is $n$.
By Lemma~\ref{He}, we know $b_{-1} = 2$, $b_{-2}=3$, \ldots, $b_{-(n-1)}=n$.
Thus we know that (\ref{e=0}) is the equality in this case.
Since $(1,-1)\not\in S$ and $(-(n-1),n-2)\not\in T$, we know 
\[
-\frac{g_0}{f_0}=-1 < \overline{u} = -\frac{u_2}{u} < -\frac{g_1}{f_1}=-\frac{n-2}{n-1} .
\]
Then (a) in (2) follows since ${\rm GCD}(u_2,u) = 1$ (Proposition~4.8 in \cite{KN}) and (\ref{det}).
Here we remark $(\gamma,\delta) \neq (1,1)$ since $u > n = f_0+f_1$.
Since $a_1 = b_{-(n-1)} = n$, (b) and (c) in (2) follow.
Since $z^u-x^{s_3}y^{t_3}$ is the negative curve, (\ref{stu}) is satisfied.
Thus $(2)_0$ is proved.

Assume that the minimal degree is bigger than $n$.
By (\ref{e=0}) and $a_2\ge n+2$,
we know 
\[
b_{-(n-1)}=n+1, \ b_{-(n-2)}=n-1, \ b_{-(n-3)}=n-2, \ \ldots, b_{-1}=2, \ a_1 = n
\]
by Lemma~\ref{He}. 
Therefore only one of $(-(n-1), 2n-1)$ and $(-(n-1),n-2)$ belongs to $T$.

Here we assume $(-(n-1), 2n-1) \in T$. 
Consider $\Gamma(-(n-1),2n-1,n-1;1,n-1,-1)$.
Let $e_i$ be the number of lattice points in $\Gamma(-(n-1),2n-1,n-1;1,n-1,-1)$
with the first component $i$.
Then we have
\[
1< e_1 < e_2 < \cdots, \ \ 1< e_{-1} < e_{-2} < \cdots
\]
There exists $(-\alpha,\beta) \in T \cap {\Bbb Z}^2$ such that $0 < \beta < \alpha$ since $\overline{u} > -1$.
We choose such $(-\alpha,\beta) \in T \cap {\Bbb Z}^2$ with $\alpha$ minimal.
By Lemma~\ref{tuduku} (1), if we sort 
$e_{-\alpha}$, $e_{-\alpha+1}$, \ldots, $e_{-1}$, $e_0$, $e_1$, \ldots, $e_{\sigma}$
($\sigma = e_{-\alpha}-\alpha-1$) into ascending order, we obtain
\[
1,2,3,\ldots, e_{-\alpha} .
\]
If $i < -\alpha$, then $b_{i} \ge e_i > e_{-\alpha}$.
If $\sigma< i  < \alpha$, then $a_i \ge e_i >e_{-\alpha}$
by the definition of $\alpha$ and Lemma~\ref{tuduku} (1).
If $i\ge  \alpha$, then $a_i>e_{-i} \ge e_{-\alpha}$.
(For $i>0$, we have
\begin{align*}
e_{-i} & = \lfloor \frac{2n-1}{n-1}i \rfloor -i+1 = \lfloor \frac{n}{n-1}i \rfloor +1 \\
a_i & \ge (n-1)i+1 .
\end{align*}
Therefore
\[
a_i-e_{-i} \ge (n-1)i-\lfloor \frac{n}{n-1}i \rfloor
\ge (n-1)i- \frac{n}{n-1}i  = \frac{(n-1)^2-n}{n-1}i > 0
\]
if $n \ge 3$.
)
Therefore, if $i > \sigma$, then $a_i > e_{-\alpha}$.
By definition, we have
\[
\Gamma(-(n-1),2n-1,n-1;1,n-1,-1) \cap [-\alpha,\sigma]_\bZ^1
\subsetneq (S\cup T)\cap [-\alpha,\sigma]_\bZ^1 
\]
since $b_{-\alpha} > e_{-\alpha}$.
Then the condition EMU for $\Delta_{\overline{t}, \overline{u}, \overline{s}}$ is satisfied.
It is a contradiction.
Therefore $(-(n-1),n-2)$ belongs to $T$.
Then we have
\[
(S\cup T) \cap [-(n-1),1]_\bZ^1
= \Gamma(-(n-1), 2n-2, n-2;1,n-1,-1)
\cap [-(n-1),1]_\bZ^1 .
\]
Thus $(3)_0$ is proved.

Next, assume that $(1)_\lambda$, $(2)_\lambda$, $(3)_\lambda$ are satisfied for some $\lambda\ge 0$.
We shall prove $(1)_{\lambda+1}$, $(2)_{\lambda+1}$, $(3)_{\lambda+1}$.
We may assume that the minimal degree of 
$\Delta_{\overline{t}, \overline{u}, \overline{s}}$ is bigger than $f_{\lambda}+f_{\lambda+1}$.
By $(3)_\lambda$, we have
\[
(S\cup T) \cap [-f_{\lambda+1},f_{\lambda}]_\bZ^1
= \Gamma(-f_{\lambda+1}, 2f_{\lambda+1}, g_{\lambda+1}; f_{\lambda},(n-1)f_{\lambda},-g_{\lambda})
\cap [-f_{\lambda+1},f_{\lambda}]_\bZ^1 .
\]
Then we know
\begin{equation}\label{3}
(S\cup T) \cap [-f_{\lambda+2},f_{\lambda+1}]_\bZ^1
\supset \Gamma(-f_{\lambda+1}, 2f_{\lambda+1}, g_{\lambda+1};f_{\lambda+1},(n-1)f_{\lambda+1},-g_{\lambda+1})
\cap [-f_{\lambda+2},f_{\lambda+1}]_\bZ^1 .
\end{equation}
By (\ref{det}), we have
\begin{align*}
& \Gamma(-f_{\lambda+1}, 2f_{\lambda+1}, g_{\lambda+1};f_{\lambda+1},(n-1)f_{\lambda+1},-g_{\lambda+1})
\cap (-f_{\lambda+2},f_{\lambda+1})_\bZ^1 \\
= &
\Gamma(-f_{\lambda+2}, 2f_{\lambda+2}, g_{\lambda+2};f_{\lambda+1},(n-1)f_{\lambda+1},-g_{\lambda+1})
\cap (-f_{\lambda+2},f_{\lambda+1})_\bZ^1 .
\end{align*}
By Lemma~\ref{tuduki} (1),
if we sort the numbers of lattice points in each columns of 
$\Gamma(-f_{\lambda+2}, 2f_{\lambda+2}, g_{\lambda+2};f_{\lambda+1},(n-1)f_{\lambda+1},-g_{\lambda+1})
\cap (-f_{\lambda+2},f_{\lambda+1})_\bZ^1$, we obtain 
$$1,2,3,\ldots, f_{\lambda+1}+f_{\lambda+2}-1 .$$
Furthermore we have
\begin{equation}\label{1}
b_{-f_{\lambda+2}} \ge f_{\lambda+1}+f_{\lambda+2}
\end{equation}
by (\ref{3}) and Lemma~\ref{tuduki} (5).
Since $(-f_{\lambda+1},g_{\lambda+1})\in T$, we have $(f_{\lambda+1},-g_{\lambda+1})\not\in S$.
Since $(-f_{\lambda+1},g_{\lambda+1}-1)\not\in T$, we have $(f_{\lambda+1},-g_{\lambda+1}+1)\in S$.
Therefore
\begin{equation}\label{2}
a_{f_{\lambda+1}} \ge (n-1)f_{\lambda+1}+g_{\lambda+1} = f_{\lambda+1}+f_{\lambda+2}
\end{equation}
by (\ref{3}) and Lemma~\ref{tuduki} (4).
Since the condition EMU is not satisfied for $\Delta_{\overline{t}, \overline{u}, \overline{s}}$,
we have
\begin{equation}\label{sono}
(S\cup T) \cap (-f_{\lambda+2},f_{\lambda+1})_\bZ^1 =
\Gamma(-f_{\lambda+2}, 2f_{\lambda+2}, g_{\lambda+2}; f_{\lambda+1},(n-1)f_{\lambda+1},-g_{\lambda+1})
\cap (-f_{\lambda+2},f_{\lambda+1})_\bZ^1 .
\end{equation}
Thus we know that the minimal degree of $\Delta_{\overline{t}, \overline{u}, \overline{s}}$ is bigger than or equal to $f_{\lambda+1}+f_{\lambda+2}$.
We have proved $(1)_{\lambda+1}$.

Assume that the minimal degree of $\Delta_{\overline{t}, \overline{u}, \overline{s}}$ is 
equal to $f_{\lambda+1}+f_{\lambda+2}$.
Then (\ref{sono}) and $a_{f_{\lambda+1}} = b_{-f_{\lambda+2}} = f_{\lambda+1}+f_{\lambda+2}$
are satisfied, and 
(\ref{3}) is the equality in this case.
Since $(-f_{\lambda+2},g_{\lambda+2}) \not\in T$ and $(f_{\lambda+1},-g_{\lambda+1}) \not\in S$, we obtain
\[
-\frac{g_{\lambda+1}}{f_{\lambda+1}} < \overline{u} = -\frac{u_2}{u} < -\frac{g_{\lambda+2}}{f_{\lambda+2}} .
\]
By ${\rm GCD}(u_2,u) = 1$ and (\ref{det}),
(a) in (2) follows.
Here we know $(\gamma,\delta)\neq (1,1)$ since $u > f_{\lambda+1}+f_{\lambda+2}$.
Since $a_{f_{\lambda+1}} = b_{-f_{\lambda+2}} = f_{\lambda+1}+f_{\lambda+2}$, (b) and (c) in (2) follow.
Since $z^u-x^{s_3}y^{t_3}$ is the negative curve, (\ref{stu}) is satisfied.
Thus $(2)_{\lambda+1}$ is proved.

Assume that the minimal degree of 
$\Delta_{\overline{t}, \overline{u}, \overline{s}}$ is bigger than $f_{\lambda+1}+f_{\lambda+2}$.
Remember that (\ref{3}), (\ref{1}), (\ref{2}) and (\ref{sono}) are satisfied in this case.
Only one of $a_{f_{\lambda+1}}$ and $b_{-f_{\lambda+2}}$ is equal to $f_{\lambda+1}+f_{\lambda+2}$.

First assume that $a_{f_{\lambda+1}}=f_{\lambda+1}+f_{\lambda+2}$.
Since $a_{f_{\lambda+1}+1}\ge f_{\lambda+1}+f_{\lambda+2}+2$ by (\ref{abplus}),
we know $b_{-f_{\lambda+2}} = f_{\lambda+1}+f_{\lambda+2}+1$.
Therefore only one of $(-f_{\lambda+2}, 2f_{\lambda+2}+1)$ and $(-f_{\lambda+2}, g_{\lambda+2})$ belongs to $T$.
Assume $(-f_{\lambda+2}, 2f_{\lambda+2}+1)\in T$.
Then we have
\[
(S\cup T)\cap [-f_{\lambda+2},f_{\lambda+1}]_\bZ^1
=
\Gamma(-f_{\lambda+2}, 2f_{\lambda+2}+1, \frac{g_{\lambda+1}f_{\lambda+2}}{f_{\lambda+1}};f_{\lambda+1}, (n-1)f_{\lambda+1}, -g_{\lambda+1})
\cap [-f_{\lambda+2},f_{\lambda+1}]_\bZ^1 .
\]
There exists $(-\alpha,\beta) \in T\cap {\Bbb Z}^2$ such that $0 < \beta < \frac{g_{\lambda+1}\alpha}{f_{\lambda+1}}$
since $u > f_{\lambda+1}+f_{\lambda+2}$.
We choose such $(-\alpha,\beta) \in T\cap {\Bbb Z}^2$ with $\alpha$ minimal.
Then we know $\alpha>f_{\lambda+2}$ and
\[
(S\cup T)\cap [-\alpha,\alpha)_\bZ^1
\supset
\Gamma(-f_{\lambda+2}, 2f_{\lambda+2}+1, \frac{g_{\lambda+1}f_{\lambda+2}}{f_{\lambda+1}};f_{\lambda+1}, (n-1)f_{\lambda+1}, -g_{\lambda+1})
\cap [-\alpha,\alpha)_\bZ^1 .
\]
Let $e$ be the number of lattice points in $\Gamma(-f_{\lambda+2}, 2f_{\lambda+2}+1, \frac{g_{\lambda+1}f_{\lambda+2}}{f_{\lambda+1}};f_{\lambda+1}, (n-1)f_{\lambda+1}, -g_{\lambda+1})$
such that the first component is $-\alpha$.
Then we know $b_i > b_{-\alpha}> e$ for $i <- \alpha$ since $(-\alpha,\beta) \in T$.
Furthermore $a_i \ge a_\alpha > e$ if $i \ge \alpha$.
(Since $\alpha>f_{\lambda+2}$,
we know $(-f_{\lambda+1}-1,g_{\lambda+1}) \not\in T$ and  $(f_{\lambda+1}+1,-g_{\lambda+1}) \in S$.
Therefore we have  $$a_\alpha \ge  (n-1)\alpha  + 2 .$$
On the other hand, since
\[
e = \lfloor \frac{2f_{\lambda+2}+1}{f_{\lambda+2}}\alpha \rfloor - \lceil \frac{g_{\lambda+1}}{f_{\lambda+1}}\alpha \rceil + 1
\le \frac{2f_{\lambda+2}+1}{f_{\lambda+2}}\alpha -  \frac{g_{\lambda+1}}{f_{\lambda+1}}\alpha  + 1 ,
\]
we know 
\[
e \le \lfloor \left( \frac{2f_{\lambda+2}+1}{f_{\lambda+2}} - \frac{g_{\lambda+1}}{f_{\lambda+1}} \right)\alpha \rfloor + 1 .
\]
Since 
\[
n-1 \ge 2 > \frac{2f_{\lambda+2}+1}{f_{\lambda+2}} - \frac{g_{\lambda+1}}{f_{\lambda+1}} ,
\]
we obtain $a_\alpha>e$ immediately.)
Since $b_{-\alpha} > e$, we know that the condition EMU for
$\Delta_{\overline{t}, \overline{u}, \overline{s}}$ is satisfied by Lemma~\ref{tuduku} (2).
It is a contradiction.
If $(-f_{\lambda+2}, g_{\lambda+2})$ belongs to $T$, $(3)_{\lambda+1}$ is satisfied.

Next we assume that $b_{-f_{\lambda+2}}=f_{\lambda+1}+f_{\lambda+2}$ and $a_{f_{\lambda+1}}\ge f_{\lambda+1}+f_{\lambda+2}+1$.
Since $(-f_{\lambda+1},g_{\lambda+1}-1)\not\in T$ and $(-f_{\lambda+1},g_{\lambda+1})\in T$ by $(3)_\lambda$, we know
$(f_{\lambda+1},-g_{\lambda+1}+1)\in S$ and $(f_{\lambda+1},-g_{\lambda+1})\not\in S$.
Then we know that $S$ contains $(f_{\lambda+1}, (n-1)f_{\lambda+1}+1)$.
We know 
\[
S \cup T \supset 
\Gamma(-f_{\lambda+2}-1,2(f_{\lambda+2}+1), \frac{g_{\lambda+1}(f_{\lambda+2}+1)}{f_{\lambda+1}};f_{\lambda+1},(n-1)f_{\lambda+1}+1,-g_{\lambda+1}) \cap [-f_{\lambda+2}-1, f_{\lambda+1}]_\bZ^1 .
\]
Let $e'_i$ be the number of the lattice points in $\Gamma(-f_{\lambda+2}-1,2(f_{\lambda+2}+1), \frac{g_{\lambda+1}(f_{\lambda+2}+1)}{f_{\lambda+1}};f_{\lambda+1},(n-1)f_{\lambda+1}+1,-g_{\lambda+1})$
with the first component $i$.
Then $e'_{f_{\lambda+1}} = f_{\lambda+1}+f_{\lambda+2}+1$,
$e'_{-f_{\lambda+2}} = f_{\lambda+1}+f_{\lambda+2}$,
$e'_{-f_{\lambda+2}-1} = f_{\lambda+1}+f_{\lambda+2}+2$.
The last equality follows from $(-f_{\lambda+2}-1, g_{\lambda+2}+1) \in \Gamma(-f_{\lambda+2}-1,2(f_{\lambda+2}+1), \frac{g_{\lambda+1}(f_{\lambda+2}+1)}{f_{\lambda+1}};f_{\lambda+1},(n-1)f_{\lambda+1}+1,-g_{\lambda+1})$.
(If not, the point $(-f_{\lambda+2}-1, g_{\lambda+2}+1)$ is in the interior of the cone 
spanned by $(-f_{\lambda+1}, g_{\lambda+1})$ and $(-f_{\lambda+2}, g_{\lambda+2})$.
It is impossible since $f_{\lambda+1} \ge f_1 = n-1 \ge 2$ and (\ref{det}).)
Here we remark 
\begin{align*}
& \Gamma(-f_{\lambda+2}-1,2(f_{\lambda+2}+1), \frac{g_{\lambda+1}(f_{\lambda+2}+1)}{f_{\lambda+1}};f_{\lambda+1},(n-1)f_{\lambda+1}+1,-g_{\lambda+1}) 
\cap (-f_{\lambda+2},f_{\lambda+1})_\bZ^1\\
= & \Gamma(-f_{\lambda+2},2f_{\lambda+2},g_{\lambda+2};
f_{\lambda+1},(n-1)f_{\lambda+1},-g_{\lambda+1})
\cap (-f_{\lambda+2},f_{\lambda+1})_\bZ^1 .
\end{align*}
Therefore, if we sort $e'_{-f_{\lambda+2}+1}, e'_{-f_{\lambda+2}+2}, \ldots, e'_{f_{\lambda+1}-1}$
into ascending order, we have $$1,2,\ldots, f_{\lambda+1}+f_{\lambda+2}-1$$
by Lemma~\ref{tuduki} (1).
Since the condition EMU for $\Delta_{\overline{t}, \overline{u}, \overline{s}}$ is not satisfied, we know
\begin{align*}
& (S \cup T) \cap [-f_{\lambda+2}-1, f_{\lambda+1}]_\bZ^1 \\
= & 
\Gamma(-f_{\lambda+2}-1,2(f_{\lambda+2}+1), \frac{g_{\lambda+1}(f_{\lambda+2}+1)}{f_{\lambda+1}};f_{\lambda+1},(n-1)f_{\lambda+1}+1,-g_{\lambda+1}) \cap [-f_{\lambda+2}-1, f_{\lambda+1}]_\bZ^1 .
\end{align*}
Here remark that $\Gamma(-f_{\lambda+2}-1,2(f_{\lambda+2}+1), \frac{g_{\lambda+1}(f_{\lambda+2}+1)}{f_{\lambda+1}};f_{\lambda+1},(n-1)f_{\lambda+1}+1,-g_{\lambda+1})$
satisfies the assumption in Lemma~\ref{tuduku}.
There exists $(-\alpha,\beta) \in T\cap {\Bbb Z}^2$ such that $0 < \beta < \frac{g_{\lambda+1}\alpha}{f_{\lambda+1}}$.
We choose such $(-\alpha,\beta) \in T\cap {\Bbb Z}^2$ with $\alpha$ minimal.
Then we know $\alpha > f_{\lambda+2}+1$ and 
\[
(S\cup T)\cap [-\alpha,\alpha)_\bZ^1
\supset
\Gamma(-f_{\lambda+2}-1,2(f_{\lambda+2}+1), \frac{g_{\lambda+1}(f_{\lambda+2}+1)}{f_{\lambda+1}};f_{\lambda+1},(n-1)f_{\lambda+1}+1,-g_{\lambda+1})
\cap [-\alpha,\alpha)_\bZ^1 .
\]
Then $b_i > b_{-\alpha} > e'_{-\alpha}$ if $i < -\alpha$.
Furthermore $a_i > e'_{-i} > e'_{-\alpha}$ if $i \ge \alpha$.
Since $b_{-\alpha} > e'_{-\alpha}$, we know that the condition EMU for
$\Delta_{\overline{t}, \overline{u}, \overline{s}}$ is satisfied by Lemma~\ref{tuduku}.
It is a contradiction.
We have completed the proof of $(3)_{\lambda+1}$.

We have completed the proof of Claim~\ref{classify} and 
$(1) \Rightarrow (2)$ in Theorem~\ref{classification}.

\vspace{2mm}

Assume (2) in Theorem~\ref{classification}.
Then we have
\[
-\frac{g_{\lambda}}{f_{\lambda}} < \overline{u} = -\frac{u_2}{u} < -\frac{g_{\lambda+1}}{f_{\lambda+1}} 
\]
and $u > f_\lambda+f_{\lambda+1}$.
By (a), (b) and (c) in (2), we know
\[
(S\cup T) \cap (-f_{\lambda+1},f_\lambda]^1_\bZ = 
\Gamma(-f_{\lambda+1},2f_{\lambda+1},g_{\lambda+1};f_{\lambda},(n-1)f_{\lambda},-g_{\lambda}) \cap (-f_{\lambda+1},f_\lambda]^1_\bZ
\]
and $-b_{f_{\lambda+1}} = f_\lambda+f_{\lambda+1}$.
By Lemma~\ref{tuduki}, if we sort $1, a_1, a_2, \ldots, a_{f_{\lambda}}, b_{-1}, b_{-2}, \ldots, b_{-f_{\lambda+1}}$ into ascending order, we obtain
\[
1,2,\ldots,  f_{\lambda}+f_{\lambda+1}-1, f_{\lambda}+f_{\lambda+1}, f_{\lambda}+f_{\lambda+1} .
\]
Thus the condition EMU is not satisfied for $\Delta_{\overline{t}, \overline{u}, \overline{s}}$.
(a) of (1) follows from (d) of (2).
(c) of (1) is clear.
We have completed the proof of  Theorem~\ref{classification}.
\qed

\begin{Remark}\label{indep}
\begin{rm}
Let $a$, $b$, $c$ be pairwise coprime positive integers such that
$\Delta_{\overline{t}, \overline{u}, \overline{s}}$ satisfies (2) in Theorem~\ref{classification} with $n$, $\lambda$, $\gamma$, $\delta$.
Let $a_i$ and $b_j$ be integers defined in (\ref{defab}).
Let $e_i$ be the number of lattice points in 
$\Gamma(-f_{\lambda+1}, 2f_{\lambda+1}, g_{\lambda+1};f_{\lambda},(n-1)f_{\lambda},-g_{\lambda})$ with the first component $i$.
Then we know
\begin{itemize}
\item
$a_i = e_i$ for $i = 1, 2, \ldots, f_\lambda$,
\item
$b_{-i} = e_{-i}$ for $i = 1, 2, \ldots, f_{\lambda+1}-1$,
\item
$a_{f_\lambda} = b_{-f_{\lambda+1}} = e_{-f_{\lambda+1}}-1 = f_\lambda+f_{\lambda+1}$.
\end{itemize}
Remark that 
\[
a_1, a_2, \ldots, a_{f_\lambda}, b_{-1}, b_{-2}, \ldots, b_{-f_{\lambda+1}}
\]
are independent of $\gamma$ and $\delta$.
\end{rm}
\end{Remark}

\begin{Remark}\label{sonzai}
\begin{rm}
For given $\lambda \in \bNo$ and $\gamma, \delta \in \bN$ with $(\gamma, \delta) \neq (1,1)$ and ${\rm GCD}(\gamma,\delta) = 1$,
it is possible to find pairwise coprime positive integers $a$, $b$, $c$ satisfying
(a), (b), (c), (d) in (2) of Theorem~\ref{classification} as follows.

We put $u = \gamma f_\lambda + \delta f_{\lambda+1}$, $u_2 = \gamma g_\lambda + \delta g_{\lambda+1}$, $\overline{u} = -u_2/u$,
$s' = \frac{(n-1)f_\lambda+1}{f_\lambda}$,
$t' = -\frac{2f_{\lambda+1}+1}{f_{\lambda+1}}$.
Consider the triangle $\Delta_{t', \overline{u}, s'}$.
Then the sequence $\ell'_1, \ell'_2, \ldots$ in Definition~\ref{Keu}
is equal to 
\[
1,2,\ldots, f_\lambda+f_{\lambda+1}-1, f_\lambda+f_{\lambda+1}+1, f_\lambda+f_{\lambda+1}+1, \ldots .
\]
Then the condition EMU for $\Delta_{t', \overline{u}, s'}$ is satisfied by Lemma~\ref{He}.
Consider the convex hull $P$ of $\Delta_{t', \overline{u}, s'} \cap \bZ^2$.
Let $B$ be the number of lattice points in the boundary of $P$.
Let $I$ be the number of lattice points in the interior of $P$.
Since the condition EMU is satisfied,
we know
\[
B+I \ge 1 + (1+2+\cdots+u)+1 = 2 + \frac{u(u+1)}{2}
\]
It is easy to see $B\le u+1$.
Then, by Pick's theorem, we know 
\[
|\Delta_{t', \overline{u}, s'}| \ge
|P| = \frac{B}{2} + I-1 \ge \frac{u+1}{2} + \left\{ 2 + \frac{u(u+1)}{2} - (u+1) \right\} -1 = \frac{u^2+1}{2} > \frac{u^2}{2} .
\]
Then we have
\[
\frac{1}{s'-\overline{u}} + \frac{1}{\overline{u}-t'}<1 .
\]
We choose sufficiently near $\overline{s} < s'$ and $\overline{t} < t'$ 
(see Corollary~\ref{abc}), 
we can find $a$, $b$, $c$ satisfying the requirement.
(The sequence $\ell'_1, \ell'_2, \ldots$ of $\Delta_{\overline{t}, \overline{u}, \overline{s}}$ is equal to 
\[
1,2,\ldots,  f_\lambda+f_{\lambda+1}-1, f_\lambda+f_{\lambda+1}, f_\lambda+f_{\lambda+1}, \ldots .
\]
Therefore the minimal degree is equal to $f_\lambda+f_{\lambda+1}$.)
\end{rm}
\end{Remark}

\section{A proof of Theorem~\ref{Keuprop}}

We shall prove Theorem~\ref{Keuprop} in this section.

\begin{Lemma}\label{21}
Let $K$ be a field of characteristic $0$.
Let $n$ and $\lambda$ be integers such that $n\ge 3$ and $\lambda\ge 0$.

Then there exists pairwise coprime positive integers $a$, $b$, $c$ satisfying
following conditions:
\begin{enumerate}
%\item
%(\ref{abcjyouken}) holds,
\item
$R_s({\frak p})$ is not Noetherian.
\item
The condition (2) in Theorem~\ref{classification} is satisfied with $n$, $\lambda$, $\gamma=2$, $\delta=1$.
\item
The sequence $\ell'_1, \ell'_2, \ldots, \ell'_u$ in Definition~\ref{Keu} is equal to
\begin{equation}\label{sequence}
1, 2, \ldots, f_{\lambda}+f_{\lambda+1}-1, f_{\lambda}+f_{\lambda+1}, f_{\lambda}+f_{\lambda+1}, 
f_{\lambda}+f_{\lambda+1}+2,f_{\lambda}+f_{\lambda+1}+3, \ldots,
2f_{\lambda}+f_{\lambda+1} 
\end{equation}
if $\lambda > 0$, and 
\[%begin{equation}\label{sequence2}
1, 2, \ldots, f_{\lambda}+f_{\lambda+1}-1, f_{\lambda}+f_{\lambda+1}, f_{\lambda}+f_{\lambda+1}
\]%end{equation}
if $\lambda = 0$.
\end{enumerate}
\end{Lemma}

\proof
First we assume $\lambda= 0$.
By Remark~\ref{sonzai}, we can find 
pairwise coprime positive integers $a$, $b$, $c$ satisfying 
(2) in Theorem~\ref{classification} with $n$, $\lambda=0$, $\gamma=2$, $\delta=1$.
The minimal degree is equal to $f_{0}+f_{1}=n$.
In this case, $u = 2f_0+f_1=n+1$.
Therefore (3) in Lemma~\ref{21} is satisfied by Lemma~\ref{He}.
In this case, by Remark~\ref{indep},
\[
a_1=n, \ b_{-1} = 2, \ b_{-2} = 3, \ \ldots, \ b_{n-1} = n.
\]
We know that 
$R_s({\frak p})$ is not Noetherian by Theorem~1.2 in \cite{GK}.

Next we assume $\lambda> 0$.
Consider $\Gamma(-f_{\lambda+1},2f_{\lambda+1},g_{\lambda+1};
f_{\lambda},(n-1)f_{\lambda},-g_{\lambda})$.
Let $e_i$ be the number of lattice points in $\Gamma(-f_{\lambda+1},2f_{\lambda+1},g_{\lambda+1};f_{\lambda},(n-1)f_{\lambda},-g_{\lambda})$ with the first component $i$.
If we sort $\{ e_i \mid i \in \bZ \}$ into ascending order,
we obtain $$1,2,3,\ldots$$ by Lemma~\ref{tuduki} (1).
Put $u' =-\frac{2g_{\lambda}+g_{\lambda+1}}{2f_{\lambda}+f_{\lambda+1}}$.
Consider $\Delta_{-2,u',(n-1)}$.
Then we have
\[
(S \cup T)\cap \bZ^2 \subset \Gamma(-f_{\lambda+1},2f_{\lambda+1},g_{\lambda+1};f_{\lambda},(n-1)f_{\lambda},-g_{\lambda}) .
\]
We define $a_i$, $b_i$ as in (\ref{defab}).

Then we have $a_i \le e_i$ and $b_{-i} \le e_{-i}$ for $i > 0$.
%We define 
%\begin{equation}\label{qq}
%q'_1, q'_2, \ldots, q'_{2f_{\lambda}+f_{\lambda+1}}
%\end{equation}
%as in (\ref{q}).
Remark
\[
\cdots > b_{-(i+1)} > b_{-i} > \cdots > b_{-2} > b_{-1} > 1 < a_1<a_2<\cdots<a_i<a_{i+1}<\cdots 
\]
and remember (\ref{seq}).
The sequence $\ell_1$, $\ell_2$, \ldots, $\ell_{2f_\lambda+f_{\lambda+1}} $ defined in Definition~\ref{Keu} 
is equal to 
\[
a_1, \ a_2 , \ \ldots, \ a_j, \ b_{-k}, \ \ldots, \ b_{-2}, \ b_{-1}, \ 1
\]
for some $j$, $k$ such that $j+k+1=2f_{\lambda}+f_{\lambda+1}$.
Since 
\[
a_{f_\lambda} = b_{-f_{\lambda+1}} = f_{\lambda}+f_{\lambda+1} ,
\]
we know $j \ge f_\lambda$ and $k \ge f_{\lambda+1}$.
We have
\[
\mbox{$a_i = e_i$ for $i = 1, 2, \ldots, 2f_{\lambda}+f_{\lambda+1}-1$}
\]
and
\[
b_{-i} = \left\{
\begin{array}{ll}
e_{-i} & i = 1, 2, \ldots, f_{\lambda+1}-1 , \\
e_{-f_{\lambda+1}}-1 & i =  f_{\lambda+1} , \\
e_{-i} & i = f_{\lambda+1}+1, \ldots, f_{\lambda}+f_{\lambda+1}-1 .
\end{array}
\right.
\]
Then
we know that the sequence 
$\ell'_1$, $\ell'_2$, \ldots, $\ell'_{2f_\lambda+f_{\lambda+1}}$ for $\Delta_{-2,u',(n-1)} $ defined in Definition~\ref{Keu} 
is equal to (\ref{sequence})
by Lemma~\ref{tuduki} (1).
Here remark that this sequence does not end at $f_{\lambda}+f_{\lambda+1}$
by $f_\lambda > 1$ (since $\lambda > 0$).
Therefore 
$\Delta_{-2,u',(n-1)}$ has a column with $2f_{\lambda}+f_{\lambda+1}$ lattice points.
Then it is easy to see that there exists 
pairwise coprime positive integers $a$, $b$, $c$ satisfying following conditions:
\begin{itemize}
\item
$\Delta_{-2,u',(n-1)} \cap \bZ^2 = \Delta_{\overline{t}, \overline{u}, \overline{s}} \cap \bZ^2$, in particular $\overline{t}<-2$, $u' = \overline{u}$, $n-1<\overline{s}$.
\item
$|\Delta_{\overline{t}, \overline{u}, \overline{s}}| > (2f_{\lambda}+f_{\lambda+1})^2/2$.
\end{itemize}
By the second condition, we know that $z^u-x^{s_3}y^{t_3}$ is a negative curve.
Thus (2) in Theorem~\ref{classification} is satisfied for $n$, $\lambda$, $\gamma=2$, $\delta=1$.
Here $u = 2f_{\lambda}+f_{\lambda+1}$, $u_2 = 2g_{\lambda}+g_{\lambda+1}$.

Assume that $R_s({\frak p})$ is Noetherian.
By Proposition~\ref{ch=0} (Theorem~1.1 in \cite{KN}),
there exists a Laurent polynomial
\[
g \in \sum_{(\alpha,\beta) \in \Delta_{\overline{t}, \overline{u}, \overline{s}} \cap \bZ^2} K v^\alpha w^\beta \cap (v-1,w-1)^u \subset K[v^{\pm1}, w^{\pm1}]
\]
such that the constant term of $g$ is not $0$.
We know that the coefficient of $v^uw^{-u_2}$ in $g$ is not $0$
since this curve does not meet the negative curve.
Then we know that $g$ is irreducible in $K[v^{\pm1}, w^{\pm1}]$ by Lemma~2.3 \cite{GAGK1}.
Consider the convex hull $P$ of $\Delta_{\overline{t}, \overline{u}, \overline{s}} \cap \bZ^2$. 
Let $Q$ be the Newton polygon of $g$.
Then $Q \subset P$.
The number of lattice points in the boundary of $P$ is $u+1$ since 
 $\Delta_{-2,u',(n-1)} \cap \bZ^2 = \Delta_{\overline{t}, \overline{u}, \overline{s}} \cap \bZ^2$.
 Here remark that each column has just $1$ point in the boundary of $P$ as in the picture below.
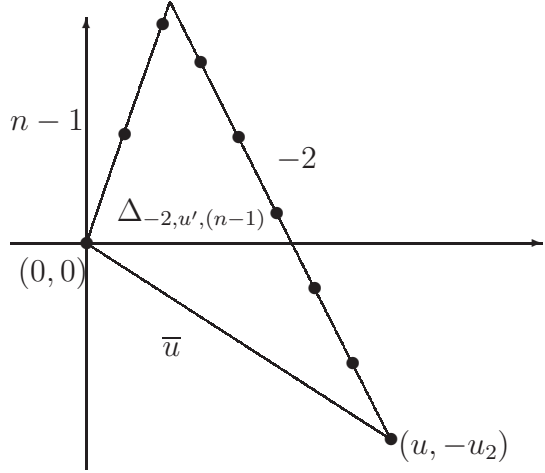
\[
{
\setlength\unitlength{1truecm}
  \begin{picture}(7,6)(-2,-3)
  \put(-1,0){\vector(1,0){7}}
  \put(0,-3){\vector(0,1){6}}
\qbezier (0,0) (2,-1.3) (4,-2.6)
\qbezier (0,0) (0.5,1.5) (1.1,3.2)
\qbezier (1.1,3.2) (2.625,0.15) (4,-2.6)
\put(-0.1,-0.1){$\bullet$}
%\put(0.15,0.625){$\bullet$}
\put(0.4,1.35){$\bullet$}
%\put(0.65,2.075){$\bullet$}
\put(0.9,2.8){$\bullet$}
\put(3.9,-2.7){$\bullet$}
%\put(3.65,-2.2){$\bullet$}
\put(3.4,-1.7){$\bullet$}
%\put(3.15,-1.2){$\bullet$}
\put(2.9,-0.7){$\bullet$}
%\put(2.65,-0.2){$\bullet$}
\put(2.4,0.3){$\bullet$}
%\put(2.15,0.8){$\bullet$}
\put(1.9,1.3){$\bullet$}
%\put(1.65,1.8){$\bullet$}
\put(1.4,2.3){$\bullet$}
%\put(1.15,2.8){$\bullet$}
\put(-0.9,-0.5){$(0,0)$}
\put(4.1,-2.8){$(u, -u_2)$}
  \put(1,-1.5){$\overline{u}$}
  \put(2.5,1){$-2$}
  \put(-1,1.5){$n-1$}
      \put(0.4,0.3){$\Delta_{-2,u',(n-1)}$}
  \end{picture}
}
\]
%Put $u = 2f_{\lambda}+f_{\lambda+1}$.
The number of lattice points in $P$ is
\[
1 + (1+2+\cdots+(2f_{\lambda}+f_{\lambda+1}))-1
=\frac{u(u+1)}{2}
\]
by (\ref{sequence}).
Therefore the number of lattice points in the interior of $P$ is
\[
\frac{u(u+1)}{2} - (u+1) = \frac{u(u-1)}{2} -1.
\]
Then, by Pick's theorem, we have
\[
|Q| \le |P| = \frac{u+1}{2} + \left( \frac{u(u-1)}{2} -1 \right) -1 = \frac{u^2-3}{2}<
\frac{u^2}{2} .
\]
Since $g$ is irreducible, $g$ is a $u$-nct in the sense of \cite{K42}.
The number $I$ of lattice points in the interior of $Q$ is less than or equal to that of $P$.
Therefore we have
\[
I \le \frac{u(u-1)}{2} -1 .
\]
It contradicts to Lemma~4.1 in \cite{K42}.
Therefore $R_s({\frak p})$ is not Noetherian.

We have completed the proof of Lemma~\ref{21}.
\qed

\vspace{2mm}

Now we start to prove Theorem~\ref{Keuprop}.
If the condition EMU is satisfied for $(a,b,c)$,
the symbolic Rees ring of ${\frak p}$ is Noetherian by Propositoin~4.6 in \cite{KN}.
Here we shall prove the converse.

Assume that the condition EMU for $\Delta_{\overline{t}, \overline{u}, \overline{s}}$ is not satisfied.
Let $n$ and $\lambda$ be integers such that $n\ge 3$ and $\lambda \ge 0$.
We may assume that $\ell_1=n$ and $\ell_{u-1}=2$ as in Remark~\ref{rem6.1}.
We shall consider pairwise coprime positive integers $a$, $b$, $c$ 
satisfying (2) in Theorem~\ref{classification} for the fixed $n$, $\lambda$ and some $\gamma$, $\delta$.
Remark that the minimal degree is $f_{\lambda} +f_{\lambda+1}$.
We put $d = f_{\lambda} +f_{\lambda+1}$.
Consider $P_A$ and $P_B$ defined just before Lemma~\ref{EMU}.
Consider 
\begin{equation}\label{Pa}
P_A \cap (\bZ \times [0, d]_\bZ)
\end{equation}
and
\begin{equation}\label{Pb}
P_B \cap (\bZ \times [0, d]_\bZ) .
\end{equation}
 (\ref{Pa}) and (\ref{Pb}) are determined by the sequences $a_1$, $a_2$, \ldots, $a_{f_{\lambda}}=d$ and
$b_{-1}$, $b_{-2}$, \ldots, $b_{-f_{\lambda+1}}=d$.
Therefore (\ref{Pa}) and (\ref{Pb}) are independent of the choice of $\gamma$ and $\delta$ (see Remark~\ref{indep}).

By Lemma~\ref{EMU} (3) and Remark~\ref{sennryaku} (1), (2),
there exist $\eta_{A,1}, \eta_{A,2} \in {A'_{d+1}}^\times$ and $\eta_{B,1}, \eta_{B,2} \in {B'_{d+1}}^\times$, such that 
\[
\mbox{$(1+x)\eta_{A,1}\eta_{B,1} = 1 + q_1x_{f_{\lambda},d}$ in ${F'_{d+1}}^\times$}
\]
and
\[
\mbox{$w \eta_{A,2}\eta_{B,2} = 1 + q_2x_{f_{\lambda},d}$ in ${F'_{d+1}}^\times$}
\]
for some integers $q_1$, $q_2$.\footnote{The integers $q_1$ and $q_2$ depend only on $n$ and $\lambda$.
They are independent of $\gamma$ and $\delta$.
Here we assume that $\gamma$ and $\delta$ are positive integers.
(We do not assume $(\gamma,\delta) \neq (1,1)$.)

If $-d < \alpha < d$, then $\lceil \alpha\overline{u} \rceil$ is independent of $\gamma$ and $\delta$ by (\ref{det}), where $\overline{u} = -\frac{\gamma g_{\lambda} + \delta g_{\lambda+1}}{\gamma f_{\lambda} + \delta f_{\lambda+1}}$.
By (\ref{Pa}) and (\ref{Pb}), we know that the rings $F'_{d+1}$, $A'_{d+1}$, $B'_{d+1}$ are 
 independent of $\gamma$ and $\delta$.}
Then, for any integers $h_1$ and $h_2$, we obtain
\begin{equation}\label{hq}
(1+x)^{h_1} w^{h_2}
\left(\eta_{A,1}^{h_1}\eta_{A,2}^{h_2}\right)
\left(\eta_{B,1}^{h_1}\eta_{B,2}^{h_2}\right)
= (1 + q_1x_{f_{\lambda},d})^{h_1} (1 + q_2x_{f_{\lambda},d})^{h_2}
= 1+(h_1q_1+h_2q_2)x_{f_{\lambda},d}
\end{equation}
in ${F'_{d+1}}^\times$.

Put $d' = g_{\lambda} +g_{\lambda+1}$.
Now consider $\Delta_{-2,-d'/d,(n-1)}$.
Considering $\Gamma(-f_{\lambda+1},2f_{\lambda+1},g_{\lambda+1};
f_{\lambda},(n-1)f_{\lambda},-g_{\lambda})$, we know that 
the sequence $\ell'_1$, $\ell'_2$, \ldots, $\ell'_u$ given in Definition~\ref{Keu}
is $1$, $2$, \ldots, $d$.
In particular, $\Delta_{-2,-d'/d,(n-1)}$ has a column that has $d$ lattice points.
Then it is easy to see that there exists 
pairwise coprime positive integers $a$, $b$, $c$ satisfying following conditions:
\begin{itemize}
\item
$\Delta_{-2,-d'/d,(n-1)} \cap \bZ^2 = \Delta_{\overline{t}, \overline{u}, \overline{s}} \cap \bZ^2$, in particular $\overline{t}<-2$, $-d'/d = \overline{u}$, $n-1<\overline{s}$.
\item
$|\Delta_{\overline{t}, \overline{u}, \overline{s}}| > d^2/2$.
\end{itemize}
By the second condition, we know that $z^u-x^{s_3}y^{t_3}$ is a negative curve.
Then the condition EMU for $\Delta_{\overline{t}, \overline{u}, \overline{s}}$ is satisfied
and $\oo_Y(abH-uE)|_{uC} \simeq \oo_{uC}$ by Remark~\ref{criteria}, Lemma~\ref{EMU} and Remark~\ref{sennryaku}.
Therefore, by Corollary~\ref{mwofuyasu}, we have
$\oo_Y(m'(abH-uE))|_{m'uC} \simeq \oo_{m'uC}$ for any $m'>0$.
Hence there exists $\eta_A \in {A'_{m'u}}^\times$ and  $\eta_B \in {B'_{m'u}}^\times$ such that
\[
\mbox{$(1+x)^{m'u}w^{-m'u_2} = \eta_A\eta_B$ in ${F'_{m'u}}^\times$}
\]
by Proposition~\ref{criterion}.
Here suppose $m' \ge 2$. 
%Then, since $m'u > d+1$, there exists $\eta_{A,1} \in {A'_{d+1}}^\times$ and  $\eta_{B,1} \in {B'_{d+1}}^\times$ such that
%\[
%\mbox{$(1+x)^{m'u}w^{-m'u_2} = \eta_{A,1}\eta_{B,1}$ in ${F'_{d+1}}^\times$} .
%\]
Then we know $m'u>d+1$ and the condition EMU is not satisfied for $m'\Delta_{\overline{t}, \overline{u}, \overline{s}}$.
Then we have 
\[
m'uq_1 + (-m'u_2)q_2 = 0
\]
by Lemma~\ref{strategy} and (\ref{hq}).
Thus we have $uq_1 -u_2q_2 = 0$.
Since ${\rm GCD}(u,u_2)= 1$, we know
\[
q_1 = qu_2 = q(g_{\lambda}+g_{\lambda+1}), \ \ q_2 = qu = q(f_{\lambda}+f_{\lambda+1})
\]
for some integer $q$.

On the other hand, consider the pairwise coprime positive integers satisfying Lemma~\ref{21}.
By the condition (3) in Lemma~\ref{21}, we know that (\ref{onlyone}) is satisfied.
Then, by Lemma~\ref{strategy2}, 
we know
\[
q_1(2f_{\lambda}+f_{\lambda+1}) + q_2(-2g_{\lambda}-g_{\lambda+1}) \neq 0 .
\]
In particular, we obtain $q\neq 0$.

Here let $\gamma$ and $\delta$ be positive integers such that
${\rm GCD}(\gamma,\delta)=1$ and $(\gamma,\delta)\neq (1,1)$.
Then we have
\begin{align*}
& q_1(\gamma f_{\lambda}+\delta f_{\lambda+1}) + q_2(-\gamma g_{\lambda}-\delta g_{\lambda+1}) \\
= & q\left\{
(g_{\lambda}+g_{\lambda+1})(\gamma f_{\lambda}+\delta f_{\lambda+1})
+(f_{\lambda}+f_{\lambda+1})(-\gamma g_{\lambda}-\delta g_{\lambda+1})
\right\} \\
= & q \times  {\rm det}
\left(
\begin{array}{cc}
\gamma f_{\lambda}+\delta f_{\lambda+1} & f_{\lambda}+f_{\lambda+1} \\
\gamma g_{\lambda}+\delta g_{\lambda+1} & g_{\lambda}+g_{\lambda+1}
\end{array}
\right) \\
= & q \times {\rm det}
\left(
\begin{array}{cc}
f_{\lambda} & f_{\lambda+1} \\
g_{\lambda} & g_{\lambda+1}
\end{array}
\right) \times
 {\rm det}
\left(
\begin{array}{cc}
\gamma & 1 \\
\delta & 1
\end{array}
\right) \\
\neq & 0
\end{align*}
by the choice of $\gamma$ and $\delta$, (\ref{det}) and $q\neq 0$.
Then, for pairwise coprime positive integers $a$, $b$, $c$ 
satisfying (2) in Theorem~\ref{classification} with $n$, $\lambda$, $\gamma$, $\delta$, 
the symbolic Rees ring $R_s({\frak p})$ is not Noetherian by Lemma~\ref{strategy}.

We have completed the proof of Theorem~\ref{Keuprop}.
\qed

\vspace{3mm}

\noindent
\begin{tabular}{l}
Taro Inagawa \\
Department of Mathematics \\
Faculty of Science and Technology \\
Meiji University \\
Higashimita 1-1-1, Tama-ku \\
Kawasaki 214-8571, Japan \\
{\tt kado@ta3.so-net.ne.jp} 
\end{tabular}

\vspace{3mm}

\noindent
\begin{tabular}{l}
Kazuhiko Kurano \\
Department of Mathematics \\
Faculty of Science and Technology \\
Meiji University \\
Higashimita 1-1-1, Tama-ku \\
Kawasaki 214-8571, Japan \\
{\tt kurano@meiji.ac.jp} \\
{\tt http://www.isc.meiji.ac.jp/\~{}kurano}
\end{tabular}


\begin{thebibliography}{99}
%\bibitem{Cox} {\sc D. A. Cox},
%{\it The Homogeneous Coordinate Ring of a Toric Variety},
%J Algebraic Geometry {\bf 4} (1995), 17--50.

\bibitem{Cowsik}  {\sc R. C. Cowsik},
{\it Symbolic powers and number of defining equations}, 
Algebra and its applications (New Delhi, 1981), 
Lecture Notes in Pure and Appl. Math. 91, Dekker, New York, 1984, 13-14.

\bibitem{C} {\sc S. D. Cutkosky}, 
{\it Symbolic algebras of monomial primes}, 
{\rm J. reine angew. Math.} {\bf 416} {\rm (1991), 71-89.}

\bibitem{CK}
{\sc S D. Cutkosky and K. Kurano}, 
{\it Asymptotic regularity of powers of ideals of points in a weighted projective plane}, 
Kyoto J. Math. {\bf 51} (2011), 25--45. 

\bibitem{Ebina}
{\sc T. Ebina},
Master theses, Meiji University 2017 (Japanese).

\bibitem{GK} {\sc J. L. Gonz\'alez and K. Karu},
{\it Some non-finitely generated Cox rings},
{\rm Compos. Math. } {\bf 152} {\rm (2016),  984--996.}

\bibitem{GAGK}
{\sc J. Gonz\'alez-Anaya, J. L. Gonz\'alez and K. Karu},
{\it Constructing non-Mori Dream Spaces from negative curves},
J, Algebra  {\bf 539}  (2019),  118--137.

\bibitem{GAGK1}
{\sc J. Gonz\'alez-Anaya, J. L. Gonz\'alez and K. Karu},
{\it Curves generating extremal rays in blowups ofweighted projective planes},
J. of London Math. Soc.  {\bf 104}  (2021),  1342--1362.


\bibitem{GAGK2}
{\sc J. Gonz\'alez-Anaya, J. L. Gonz\'alez and K. Karu},
{\it The geography of negative curves}, 
arXiv:2104.03950, 2021.

\bibitem{GAGK3}
{\sc J. Gonz\'alez-Anaya, J. L. Gonz\'alez and K. Karu},
{\it Non-existence of negative curves},
arXiv:2110.13333, 2021.

%\bibitem{G} {\sc S. Goto},
%{\it The Cohen-Macaulay symbolic Rees algebras for curve singularities},
%Memoirs of Amer. Math. Soc. {\bf 526} (1994), 1--68.
%
%\bibitem{GNSnagoya}{\sc S. Goto, K. Nishida and Y. Shimoda},
%{\it Topics on symbolic Rees algebras for space monomial curves},
%Nagoya Math. J. {\bf 124} (1991), 99--132.

%\bibitem{GHNV}{\sc S. Goto, M. Herrmann, K. Nishida and O. Villamayor},
%{\it On the structure of Noetherian symbolic Rees algebras},
%Manuscripta Math. {\bf 67} (1990), 197--225.

\bibitem{GNW} {\sc S. Goto, K. Nishida and K.-i. Watanabe},
{\it Non-{C}ohen-{M}acaulay symbolic blow-ups 
for space monomial curves and counterexamples to {C}owsik's
question}, 
{\rm Proc. Amer. Math. Soc.} {\bf 120} {\rm (1994), 383--392.}

\bibitem{H} {\sc R. Hartshorne}, 
{\it Algebraic Geometry}, 
GTM 52, Springer Verlag, 1977.

\bibitem{He} {\sc Z. He}, New examples and non-examples of Mori dream spaces when blowing up toric surfaces; arXiv:1703.00819, 2017.

\bibitem{Her} {\sc J. Herzog},
{\it Generators and relations of Abelian semigroups and semigroup rings},
Manuscripta Math. {\bf 3} (1970), 175--193.

%\bibitem{Hibi} {\sc T. Hibi},
%{\it Algebraic Combinatorics on Convex Polytopes}, 
%Carslaw Publications, Glebe NSW, Australia, 1992.

\bibitem{Hu} {\sc C. Huneke},
{\it Hilbert functions and symbolic powers},
Michigan Math. J. {\bf 34} (1987), 293--318.

%\bibitem{I1} {\sc A. Ito},
%{\it Examples of Mori dream spaces with Picard number two},
%Manuscripta Math. {\bf 145} (2014), 243--254.

%\bibitem{I2} {\sc A. Ito},
%{\it Examples of Mori dream spaces with Picard number two},
%Proceeding of 26th Seminar on commutative algebra (Japanese).

%\bibitem{Nishida} {\sc K. Kai and K. Nishida},
%{\it Finitely generated symbolic Rees rings of ideals defining certain finite sets of points in ${\bP}^2$},
%preprint.
%
%\bibitem{K8} {\sc K. Kurano},
%{\it Positive characteristic finite generation of symbolic Rees algebras and Roberts' counterexamples to the fourteenth problem of Hilbert}, 
%Tokyo J. Math. {\bf 16} (1993), 473--496.

\bibitem{K42} {\sc K. Kurano},
{\it Equations of negative curves of blow-ups of Ehrhart rings of rational convex polygons}, 
J. Algebra {\bf 590} (2022), 413-438.

\bibitem{KM} {\sc K. Kurano and N. Matsuoka},
{\it On finite generation of symbolic Rees rings of space monomial curves and existence of negative curves}, 
J. Algebra {\bf 322} (2009), 3268-3290.

\bibitem{KN} {\sc K. Kurano and K. Nishida},
{\it Infinitely generated symbolic Rees rings of space monomial curves having negative curves}, 
Michigan Math. J. {\bf 68} (2019), 405--445.

\bibitem{Matsu}
{\sc M. Matsuura},
Master theses, Meiji University 2019 (Japanese).

%
%\bibitem{Mora} {\sc M. Morales},
%{\it Noetherian symbolic blow-ups},
%J. Algebra {\bf 140} (1991), 12--25.
%
%\bibitem{MG} {\sc M. Morimoto and S. Goto},
%{\it Non-Cohen-Macaulay symbolic blow-ups for space monomial curves},
%Proc. Amer. Math. Soc.  {\bf 116} (1992), 305--311.
%
%\bibitem{Okawa} {\sc S. Okawa},
%{\it On images of Mori dream spaces}, 
%Math. Ann. {\bf 364} (2016), 1315–1342. 
%
%\bibitem{Ro} {\sc P. Roberts},
%{\it An infinitely generated Symbolic blow-ups in a power series ring and a new counterexample to Hilbert's fourteenth problem}, 
%J. Alg., {\bf 132} (1990), 461--473.
%
%\bibitem{Srinivasan} {\sc H. Srinivasan},
%{\it On finite generation of symbolic algebras of monomial primes},
%Comm. in Alg. {\bf 9} (1991), 2557--2564.
%
%\bibitem{TVAV} {\sc D. Testa, A. Varilly-Alvarado and M. Velasco}, 
%{\it Big rational surfaces}, Math. Ann. {\bf 351} (2011),
%95--107.

\bibitem{Uchi}
{\sc K. Uchisawa},
Master theses, Meiji University 2017 (Japanese).

\end{thebibliography}
\end{document}